\documentclass[11pt]{article}
\usepackage[margin=1.1in]{geometry}

\usepackage{amsmath,amssymb,amsthm,mathtools}
\usepackage{bm,bbm}
\usepackage{graphicx}
\usepackage{enumitem}
\usepackage{caption,subcaption}
\usepackage{booktabs,multirow,array,makecell}
\usepackage{xcolor}
\usepackage{tikz}
\definecolor{dtime}{RGB}{31,119,180}
\definecolor{dvalid}{RGB}{44,160,44}
\definecolor{dtree}{RGB}{214,120,20}
\definecolor{dhw}{RGB}{128,80,160}
\usetikzlibrary{arrows.meta,positioning,shapes.geometric,calc,fit,backgrounds}
\usepackage{algorithm}
\usepackage{algpseudocode}

\usepackage[numbers,sort&compress]{natbib}
\usepackage[colorlinks=true,linkcolor=blue!60!black,citecolor=green!40!black,urlcolor=blue!60!black,hypertexnames=false]{hyperref}
\hypersetup{pdftitle={QUADS: An SDP-Based Branch-and-Bound Solver for the Quadratic Assignment Problem},pdfauthor={Di Hou, Thai P.D. Nguyen, Kim-Chuan Toh}}
\numberwithin{equation}{section}

\newcommand{\algphase}[1]{\Statex\hspace{\algorithmicindent}{\color{orange!85!black}\textbf{[\,#1\,]}}}
\newcommand{\R}{\mathbb{R}}
\renewcommand{\S}{\mathbb{S}}

\newcommand{\cD}{\mathcal{D}}
\newcommand{\cE}{\mathcal{E}}
\newcommand{\cF}{\mathcal{F}}
\newcommand{\cG}{\mathcal{G}}
\newcommand{\cK}{\mathcal{K}}
\newcommand{\cO}{\mathcal{O}}
\newcommand{\cP}{\mathcal{P}}
\newcommand{\cY}{\mathcal{Y}}
\newcommand{\inprod}[2]{\langle #1,#2\rangle}
\DeclareMathOperator{\diag}{diag}
\DeclareMathOperator{\Diag}{Diag}
\DeclareMathOperator{\tr}{tr}
\DeclareMathOperator{\sym}{sym}
\DeclareMathOperator{\GL}{GL}
\DeclareMathOperator{\LAP}{LAP}
\DeclareMathOperator{\Proj}{Proj}
\newcommand{\quads}{QUADS}

\newenvironment{keywords}{\par\smallskip\noindent\textbf{Keywords.} }{\par\smallskip}
\newenvironment{MSCcodes}{\par\smallskip\noindent\textbf{MSC codes.} }{\par\smallskip}

\theoremstyle{plain}
\newtheorem{theorem}{Theorem}[section]
\newtheorem{lemma}[theorem]{Lemma}
\newtheorem{proposition}[theorem]{Proposition}
\newtheorem{corollary}[theorem]{Corollary}
\theoremstyle{definition}
\newtheorem{definition}[theorem]{Definition}
\theoremstyle{remark}
\newtheorem{remark}[theorem]{Remark}

\title{\quads: An SDP-Based Branch-and-Bound Solver for the Quadratic Assignment Problem}
\author{
Di Hou\thanks{Department of Mathematics, National
University of Singapore, Singapore 119076 ({\tt dihou@u.nus.edu}).}
\and
Thai P.D. Nguyen\thanks{Department of Mathematics, National
University of Singapore, Singapore 119076 ({\tt thai.npd@u.nus.edu}).}
\and
Kim-Chuan Toh\thanks{Department of Mathematics, National
University of Singapore, Singapore 119076 ({\tt mattohkc@nus.edu.sg}).}
}
\date{}

\begin{document}
\maketitle

\begin{abstract}
The quadratic assignment problem (QAP) is one of the most difficult problems in combinatorial optimization. It is
NP-hard even to approximate, and some QAPLIB instances with only $35$ facilities remain unsolved. Doubly nonnegative (DNN)
relaxations give some of the strongest lower bounds for the QAP. However, for $n$ facilities, they have $\cO(n^4)$ variables and
constraints, which makes them expensive at every node of a branch-and-bound tree. We propose \quads, a branch-and-bound
solver for the QAP based on DNN relaxations. To reduce the cost of a node, we solve the DNN relaxation by an
alternating direction method of multipliers (ADMM) warm started from the parent node. We derive a lower bound that is
valid at every iterate of ADMM, so that ADMM can stop as soon as the node is pruned or its bound stalls. We also
compute the more expensive bounds at a node only if they are likely to prune it, as estimated from their increases
at the nodes of the same depth processed before. To reduce the number of nodes, we propose
dual fixing by the multiplier of the DNN relaxation, which excludes assignments from the whole subtree of a node without
solving its children. We further strengthen the bounds by aggregated reformulation-linearization technique (RLT) cuts at the root. We also
exploit the structure of the data, including facilities with zero flow, symmetries, and groups of facilities with constant flow. Finally, we
search the tree in parallel on CPU cores and GPUs. Numerical experiments on QAPLIB show that, on one core, \quads\
solves five instances with $25\leq n\leq 32$ in less than two minutes each, whereas toulbar2 and Gurobi solve none of them in one
hour. With $16$ cores and two GPUs, it solves nug30 in about one minute and tai40b in less than five minutes, whereas the
shortest published times, obtained on $1728$ cores, are close to one hour and more than a day. On every QAPLIB instance that it solves and whose shortest published time exceeds five minutes, it is more than ten times faster than that time.
\end{abstract}



\section{Introduction}\label{sec:intro}

\subsection{The quadratic assignment problem}

Let $A,B\in\R^{n\times n}$ be a flow matrix and a distance matrix, where $A_{ik}$ is the flow between facilities $i$
and $k$ and $B_{jl}$ is the distance between locations $j$ and $l$. The Koopmans--Beckmann quadratic assignment problem
(QAP)~\cite{koopmans1957assignment} places each facility at a different location so that the total cost of flow times
distance is minimized, which can be formulated as
\begin{equation}\label{eq:qap}
v^{\rm QAP}\coloneqq\min\Big\{\inprod{A}{XBX^\top}\;\Big|\;X\in\Pi_n\Big\},\tag{QAP}
\end{equation}
where $\Pi_n\coloneqq\{X\in\{0,1\}^{n\times n}\mid Xe=e,\ X^\top e=e\}$ is the set of $n\times n$ permutation
matrices, $e\in\R^n$ is the vector of all ones, and $X_{ij}=1$ means that facility $i$ is placed at location $j$.
Here $\inprod{\cdot}{\cdot}$ is the trace inner product, so that the objective equals $\sum_{i,j,k,l}A_{ik}B_{jl}X_{ij}X_{kl}$.

The QAP is one of the fundamental problems of combinatorial optimization. It models facility
layout~\cite{koopmans1957assignment}, backboard wiring~\cite{steinberg1961backboard}, hospital planning, keyboard
design, and many other applications~\cite{loiola2007survey,burkard2012assignment}. Moreover, it contains the
traveling salesman, graph partitioning, and maximum clique problems as special cases~\cite{burkard2012assignment}.

The QAP is also regarded as one of the most difficult problems in combinatorial
optimization~\cite{loiola2007survey}. It is NP-hard, and it remains NP-hard to
approximate within any constant factor~\cite{sahni1976pcomplete}. This difficulty also shows in
computation. For example, the QAPLIB~\cite{burkard1997qaplib} instance nug30, with $n=30$, was first solved in
2000 on a computational grid with the equivalent of about seven years of computation on a single
workstation~\cite{anstreicher2002solving}. Moreover, as of November 2026, twenty-eight QAPLIB instances have not
been solved to proven optimality, the smallest being tai35a with $n=35$~\cite{mittelmann2026qaplib}.
Consequently, instances with about $30$ facilities have served as the benchmark of exact methods for the QAP
for more than two decades.

\subsection{SDP-based branch and bound for the QAP and its challenges}\label{sec:challenges}

Exact methods for~\eqref{eq:qap} are mostly of the branch-and-bound type. They split the set of permutations into nodes by fixing assignments and discard a node whose lower bound exceeds the objective value of the best known permutation. Their efficiency depends on the trade-off between the strength of the lower bound at a node and its cost. Bounds from linear assignment problems are inexpensive but lead to large search trees. These include the Gilmore--Lawler bound~\cite{gilmore1962optimal,lawler1963quadratic} (GLB) and the dual-ascent bounds of the reformulation-linearization technique (RLT)~\cite{adams1994improved,hahn1998lower,adams2007level2}. The convex quadratic programming bound (QPB) of Anstreicher and Brixius~\cite{anstreicher2001new} uses a convex quadratic relaxation over the doubly stochastic matrices, and it was the bound of the first solution of nug30~\cite{anstreicher2002solving}.

The doubly nonnegative (DNN) relaxation~\cite{zhao1998semidefinite,povh2009copositive} is a semidefinite programming (SDP) relaxation that lifts the assignment vector to a positive semidefinite and entrywise nonnegative matrix, and it gives some of the strongest known lower bounds. Interior-point methods are not suited to this relaxation in a branch-and-bound method. Their iterations are too expensive for a relaxation with $\cO(n^4)$ variables and constraints, and they are hard to warm start from the solution of the parent, which lies close to the boundary of the cone.

First-order splitting methods have made the DNN bound computable for QAPLIB sizes, using a facially reduced relaxation~\cite{oliveira2018admm,graham2022restricted}. Liao~\cite{liao2016branch} uses such an ADMM with early termination to compute the bounds of a parallel branch-and-bound method. Hu and Xu~\cite{hu2026exact} use its multipliers to build an exact convex quadratic reformulation, which is solved by a mixed-integer solver. Low-rank methods, such as RiNNAL and RiNNAL$+$~\cite{hou2025rinnal,hou2025rinnalplus}, are efficient when the relaxation has an optimal solution of low rank. However, the DNN relaxation of the QAP may have optimal solutions of high rank. Fujii et al.~\cite{fujii2021solving} solved several open QAPLIB instances by a parallel branch-and-bound method with the Lagrangian-DNN bound~\cite{kim2016lagrangian} on many CPU cores. These developments lead to the question
addressed in this paper:
\begin{center}
\emph{How can large QAP instances be solved exactly and efficiently with DNN bounds?}
\end{center}
Answering this question requires reducing both the time per node and the number of nodes of a branch-and-bound tree, and using parallel hardware
efficiently. These goals face four difficulties.
\begin{itemize}[leftmargin=*,itemsep=2pt]
\item \textbf{Cost of the bound at a node.} Even with a first-order method, one DNN computation takes many iterations. Each iteration requires an eigenvalue decomposition of a matrix of order $\cO(n^2)$ and a projection onto $\cO(n^4)$ polyhedral constraints, and a branch-and-bound method repeats this computation at thousands of nodes.
\item \textbf{Valid bounds from inexact solutions.} The dual objective value is a lower bound only at a dual feasible solution. Before convergence, however, the dual estimate of a first-order method is in general infeasible. Moreover, a node should stop as soon as its bound exceeds the objective value of the best known permutation, or as soon as it is clear that it cannot be pruned, and neither moment is known in advance.
\item \textbf{Size of the tree.} The total time is the time per node multiplied by the number of nodes. Even with DNN bounds, the search trees of hard QAP instances can have hundreds of thousands of nodes. Their size depends strongly on the relative gap of the bound at the root and on the facility or location branched on.
\item \textbf{Using parallel hardware.} The nodes of an SDP-based tree are expensive, and the few nodes near the root are the most expensive. Hence processing different nodes on different cores leaves most cores idle at the top of the tree. GPUs suit the many small nodes deeper in the tree. A further difficulty is the choice of the arithmetic precision. GPUs are most efficient in single precision, but the rounding errors of single precision can make a lower bound invalid.
\end{itemize}

\subsection{Overview of \quads}

\begin{figure}[!htbp]
\centering
\begin{tikzpicture}[
  font=\small,
  box/.style={rectangle,rounded corners,draw=black!60,fill=gray!8,align=center,inner sep=2pt,minimum height=14mm,text width=24.5mm},
  ctime/.style={box,draw=dtime,fill=dtime!12,very thick},
  cvalid/.style={box,draw=dvalid,fill=dvalid!12,very thick},
  ctree/.style={box,draw=dtree,fill=dtree!12,very thick},
  chw/.style={box,draw=dhw,fill=dhw!10,very thick},
  deci/.style={diamond,aspect=2.2,draw=black!60,fill=white,align=center,inner sep=1pt,font=\footnotesize},
  flow/.style={-{Latex[length=2mm]},thick},
  lab/.style={font=\scriptsize},
]
\node[box] (pool) at (31mm,0) {Open nodes};
\node[box] (cheap) at (62mm,0) {GLB and QPB\\{\scriptsize Sec.~\ref{sec:dnn}}};
\node[ctime] (admm) at (93mm,0) {DNN bound\\by ADMM\\{\scriptsize Sec.~\ref{sec:admm}}};
\node[cvalid] (cert) at (124mm,0) {Valid lower bound\\{\scriptsize Sec.~\ref{sec:cert}}};
\node[deci] (prune) at (124mm,-18mm) {Prune?};
\node[cvalid] (stop) at (93mm,-18mm) {Stop?\\{\scriptsize Sec.~\ref{sec:stop}}};
\node[ctree] (fix) at (62mm,-18mm) {Dual fixing\\{\scriptsize Sec.~\ref{sec:dualfix}}};
\node[ctree] (branch) at (31mm,-18mm) {Branch\\with orbits\\{\scriptsize Secs.~\ref{sec:exact}, \ref{sec:branch}}};
\node[ctime] (warm) at (0mm,-18mm) {Warm start\\of children\\{\scriptsize Sec.~\ref{sec:warm}}};
\node[lab,below=2.5mm of prune] (disc) {Pruned};
\draw[flow] (cheap) -- (admm);
\draw[flow] (pool) -- (cheap);
\draw[flow] (admm) -- (cert);
\draw[flow] (cert) -- (prune);
\draw[flow] (prune) -- node[right,lab]{Yes} (disc);
\draw[flow] (prune) -- node[above,lab]{No} (stop);
\draw[flow] (stop) -- node[right,lab]{No} (admm);
\draw[flow] (stop) -- node[above,lab]{Yes} (fix);
\draw[flow] (fix) -- (branch);
\draw[flow] (branch) -- (warm);
\draw[flow] (warm) |- (pool);
\begin{scope}[on background layer]
\node[fit=(pool)(cheap)(admm)(cert)(prune)(stop)(fix)(branch)(warm)(disc),inner sep=3mm] (frame) {};
\draw[dhw!60,rounded corners,fill=white,thick] ($(frame.south west)+(1.6mm,-1.6mm)$) rectangle ($(frame.north east)+(1.6mm,-1.6mm)$);
\draw[dhw,rounded corners,fill=white,very thick] (frame.south west) rectangle (frame.north east);
\end{scope}
\node[chw,text width=55mm] (hw) at (108.5mm,23mm) {Parallel workers on CPU cores and GPUs\\{\scriptsize Sec.~\ref{sec:parallel}}};
\draw[dhw,very thick] (hw.south) -- (hw.south |- frame.north);
\node[ctree] (r1) at (0mm,23mm) {Automorphisms\\{\scriptsize Sec.~\ref{sec:exact}}};
\node[ctree] (r2) at (31mm,23mm) {Clique\\bound\\{\scriptsize Sec.~\ref{sec:clique}}};
\node[ctree] (r3) at (62mm,23mm) {RLT cuts\\{\scriptsize Sec.~\ref{sec:cuts}}};
\node[fit=(r1)(r2)(r3),draw=dtree,dashed,rounded corners,inner sep=1.5mm] (root) {};
\node[lab,anchor=south west] at (root.north west) {At the root};
\draw[flow] (pool.north |- root.south) -- (pool.north);
\node[anchor=north,font=\scriptsize,align=center] at ($(frame.south)+(0,-4mm)$) {
  \textcolor{dtime}{$\blacksquare$}~time per node\quad
  \textcolor{dvalid}{$\blacksquare$}~valid bounds from inexact solutions\quad
  \textcolor{dtree}{$\blacksquare$}~size of the tree\quad
  \textcolor{dhw}{$\blacksquare$}~parallel hardware};
\end{tikzpicture}
\caption{Processing of a node in \quads.}
\label{fig:overview}
\end{figure}
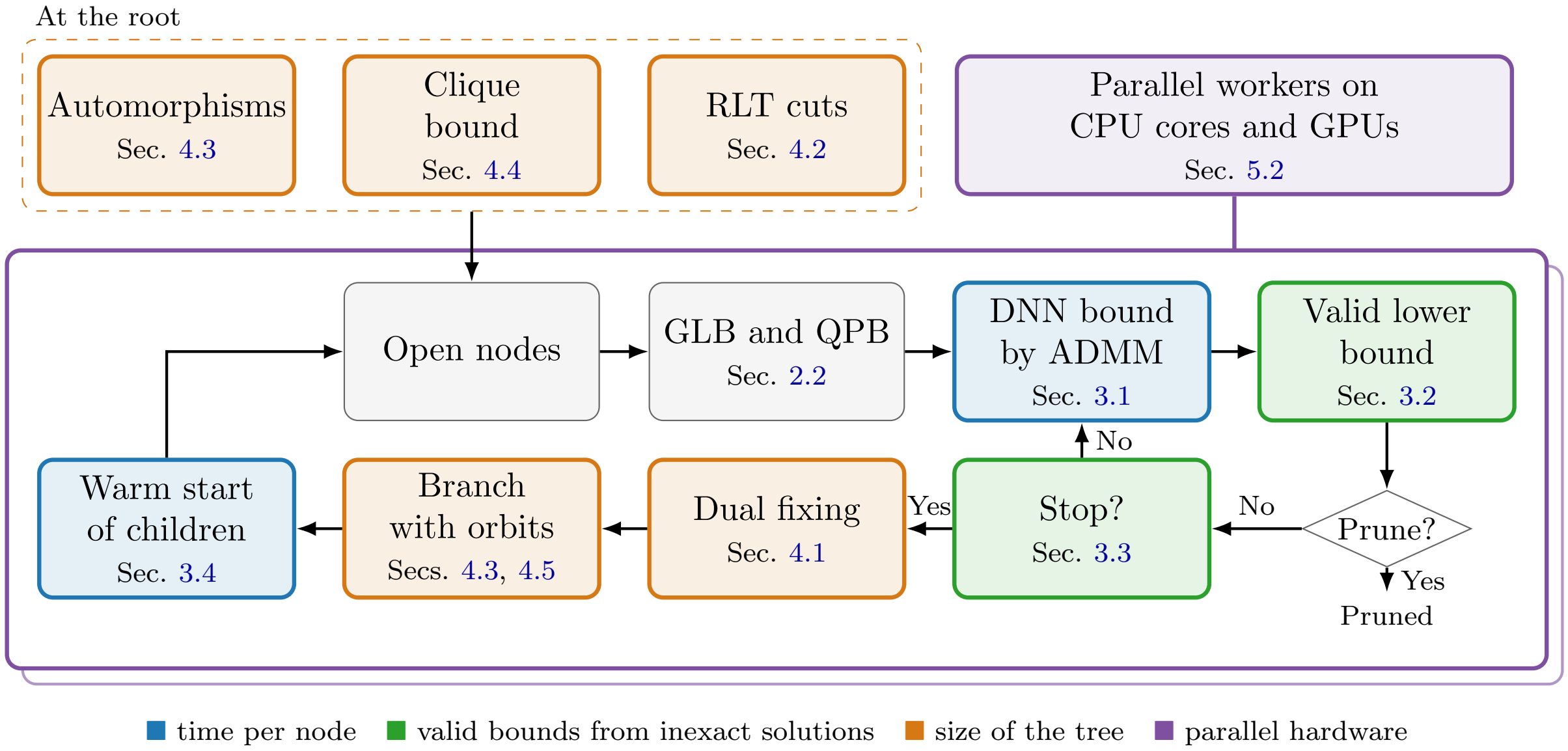

In this paper, we propose \quads\ (\textbf{QU}adratic \textbf{A}ssignment by \textbf{D}NN and \textbf{S}plitting), a
branch-and-bound solver for~\eqref{eq:qap} that addresses these difficulties. Figure~\ref{fig:overview} shows how \quads\
processes a node. Its colors mark the components developed in this paper and the difficulty of
Section~\ref{sec:challenges} that each of them addresses. At the root, \quads\ computes the automorphisms of the data,
tests the clique bound of groups of facilities with constant flow, computes the level-1 RLT bound if the GLB is close to the objective value of the best known permutation, and adds aggregated RLT cuts. Several workers on CPU cores and GPUs then search
the tree best first. At a node, \quads\ evaluates the GLB, the QPB, and the DNN bound by ADMM warm started from the
parent. The node is pruned as soon as a valid lower bound exceeds the pruning threshold. Otherwise, dual fixing
excludes assignments from its subtree. The node is then branched, and the automorphisms remove children with the same
optimal value.

\subsection{Contributions}

Our main contribution is \quads, an efficient branch-and-bound solver for the QAP based on DNN relaxations. It
makes the DNN bound affordable at every node of a branch-and-bound tree. Its components and our numerical results give
the following contributions.
\begin{enumerate}[leftmargin=*,itemsep=2pt]
\item We reduce the SDP-RLT relaxation of the QAP to a DNN relaxation of order $n^2+1$ without slack variables. We solve it by a two-block ADMM whose projections have closed forms, and every child is warm started from the final variables of its parent (Sections~\ref{sec:dnn} and~\ref{sec:bounds}).
\item We derive a valid lower bound from every multiplier of ADMM. This bound is a lower estimate of the dual function of a relaxation whose polyhedral part is the level-1 RLT without its symmetry constraints. Its supremum equals the optimal value of the DNN relaxation. Hence ADMM can stop as soon as the bound prunes the node or a forecast predicts that the node will not be pruned (Section~\ref{sec:bounds}).
\item We propose dual fixing by the multiplier of the DNN relaxation. It bounds every child of a node by a linear assignment problem with one fixed assignment, before the child is solved. Every assignment whose child bound exceeds the pruning threshold is excluded from the whole subtree of the node (Section~\ref{sec:reduce}).
\item We strengthen the bounds by summing the level-2 RLT over groups of facilities, so that its cubic variables cancel and only entries of the lifted matrix remain. After the root, their multipliers are frozen in the cost matrix, so that no node separates cuts or solves a larger relaxation (Section~\ref{sec:reduce}).
\item As a minor contribution, we incorporate various techniques in the literature to exploit the structure of the data, such as facilities with zero flow, orbital branching, and the pruning of duplicate and dominated nodes. If the free facilities form a few classes of interchangeable facilities, we impose the semidefinite constraint only on a block whose order is proportional to the number of classes. The clique bound proves optimality at the root on several QAPLIB instances (Section~\ref{sec:reduce}).
\item \quads\ runs ADMM in single precision on CPU cores and GPUs and evaluates every bound in double precision, so that every bound remains valid. A GPU processes many nodes at once, and idle cores accelerate the few expensive nodes at the top of the tree (Section~\ref{sec:impl}).
\item We conduct extensive numerical experiments on QAPLIB. On one core, \quads\ solves five instances with $25$ to $32$
  facilities in less than two minutes each, whereas toulbar2 and Gurobi solve none of them in one hour. Using
  $16$ cores and two GPUs, it solves nug30 in $67$ seconds and tai40b in $272$ seconds, whereas the shortest published
  times, obtained on $1728$ cores, are $3140$ and $105000$ seconds. On every QAPLIB instance that it solves and whose shortest published time exceeds five minutes, it is more than ten times faster than that time (Section~\ref{sec:num}).
\end{enumerate}

Section~\ref{sec:framework} describes the branch-and-bound framework, and Section~\ref{sec:conclusion} concludes the paper.

\subsection{Notation}\label{sec:notation}

The sets of real numbers and of integers are $\R$ and $\mathbb{Z}$. $\S^N$ is the space of real symmetric $N\times N$ matrices, and $\S^N_+$ is its cone of positive
semidefinite matrices. For matrices of the same size, $\inprod{X}{Y}\coloneqq\tr(X^\top Y)$ is the trace inner product and $\|\cdot\|$ is the
Frobenius norm. Inequalities between matrices, such as $Y\ge0$, hold entrywise. We write $e$ for the all-ones vector of the appropriate dimension and
$[m]\coloneqq\{1,\dots,m\}$. For a square matrix $M$, $\diag(M)$ is its diagonal and $\sym(M)\coloneqq(M+M^\top)/2$. For a vector $v$, $\Diag(v)$ is the
diagonal matrix with diagonal $v$. For $M\in\S^N$, $\lambda_{\max}(M)$ is the largest eigenvalue. For a closed convex set $\mathcal K$, $\Proj_{\mathcal K}$ is
the Euclidean projection onto $\mathcal K$.

Let $m$ be a positive integer. An assignment $p=(i,j)\in\cE\coloneqq[m]\times[m]$ places facility $i$ at location $j$. For $X\in\R^{m\times m}$, the vector $x\in\R^{m^2}$ stacks the columns
of $X$, and its entries are indexed by the assignments, so that $x_p=X_p\coloneqq X_{ij}$. Rows and columns of a matrix in $\S^{m^2+1}$ are indexed by $0$
and by the assignments. We identify a permutation $\pi$ of $[m]$ with the matrix $X\in\Pi_m$ with $X_{i\pi(i)}=1$.
For $H\in\R^{m\times m}$, the optimal value of the linear assignment problem with cost matrix $H$ is
\begin{equation}\label{eq:lap}
\LAP(H)\coloneqq\min\big\{\inprod{H}{X}\;\big|\;X\in\Pi_m\big\}=\max\big\{e^\top u+e^\top v\;\big|\;u,v\in\R^m,\ ue^\top+ev^\top\le H\big\},
\end{equation}
where the equality follows from linear programming duality, since the doubly stochastic matrices form the convex hull of $\Pi_m$~(Birkhoff--von Neumann). The problem is solved
in $\cO(m^3)$ operations by the Hungarian method~\cite{kuhn1955hungarian}. Every feasible solution $(u,v)$ of the dual problem gives the lower bound
$e^\top u+e^\top v$ on $\LAP(H)$, and $H-ue^\top-ev^\top\ge0$ is the matrix of its reduced costs.

\section{The branch-and-bound framework}\label{sec:framework}

We present the branch-and-bound framework of \quads\ and the DNN relaxation that gives its lower
bounds. Sections~\ref{sec:bounds}--\ref{sec:impl} describe the components that \quads\ adds to this framework.

\subsection{Node subproblem}\label{sec:node}

A node of the search tree is the set of the permutations $\pi$ with $\pi(i)=j$ for every assignment $(i,j)$ fixed at the node and $\pi(i)\neq j$ for every assignment $(i,j)$ excluded at the node. The root fixes and excludes no assignment. A node inherits the fixed and the excluded assignments of its parent. The depth of a node is the number of its fixed assignments. The facilities and the locations that belong to no fixed assignment are called free, and $m$ denotes their number. With a slight abuse of notation, let $A,B\in\R^{m\times m}$ be the flow and distance
submatrices of the free facilities and locations. Let $C\in\R^{m\times m}$ be the matrix whose entry $C_{ij}$ is the cost of the flow between the fixed facilities at their locations and the free facility $i$ at the free location $j$. Let $c_0$ be the cost among the fixed assignments. The assignments that are not excluded at the node are called allowed, and $\cD\subseteq\cE$ denotes their set. The node subproblem is
\begin{equation}\label{eq:node}
v^{\rm node}\coloneqq\min\Big\{\,\inprod{A}{XBX^\top}+\inprod{C}{X}+c_0\;\Big|\;X\in\Pi_m,\ X_{ij}=0\ \ \forall\,(i,j)\notin\cD\Big\}.
\end{equation}
Exchanging $A$ with $B$, $C$ with $C^\top$, and $\cD$ with $\{(j,i):(i,j)\in\cD\}$ gives an equivalent problem in $X^\top$, in which the facilities and the locations exchange their roles.

\subsection{Lower bounds}\label{sec:dnn}

\quads\ uses four lower bounds. It computes the GLB at every node, the level-1 RLT bound only at the root, and the QPB and the DNN bound at the nodes chosen by the rule at the end of this section. We first write~\eqref{eq:node} in a lifted form, in which these bounds are stated.

\paragraph{Lifted form.} For a permutation matrix $X$, one has $\inprod{A}{XBX^\top}=x^\top(B\otimes A)x$ and
$\inprod{C}{X}=c^\top x$, where $c$ stacks the columns of $C$. Since $x_p^2=x_p$ for a permutation matrix, the objective
of~\eqref{eq:node} equals $c_0+\inprod{L}{yy^\top}$ with $y=(1;x)\in\R^{N}$, $N=m^2+1$, and
\begin{equation}\label{eq:L}
L\coloneqq\begin{pmatrix}0 & 0\\ 0 & \sym(B\otimes A)+\Diag(c)\end{pmatrix}\in\S^N.
\end{equation}
The vector $y=(1;x)$ is called the lift of $X$, and $L$ is called the lifted cost matrix.
Let $\mathcal{A}\coloneqq(I_m\otimes e^\top;\,e^\top\otimes I_m)\in\R^{2m\times m^2}$, so that
$\mathcal{A}x=e$ states that every row sum and every column sum of $X$ equals one. Hence~\eqref{eq:node} is the binary
quadratic program
\begin{equation}\label{eq:Y}
v^{\rm node}=\min\big\{y^\top Ly+c_0\;\big|\;y\in\cY\big\},\qquad
\cY\coloneqq\big\{(1;x)\;\big|\;\mathcal{A}x=e,\ x\in\{0,1\}^{m^2},\ x_p=0\ \ \forall\,p\notin\cD\big\},
\end{equation}
where $\cY$ is the set of the lifts of the permutations of the node. The constraints $x_p=0$ for the excluded
assignments $p\notin\cD$ allow the entry $p$ of $y$ and the row and the column $p$ of $L$ and of the lifted matrices
below to be deleted. For simplicity, we state the
relaxations and bounds of this section without this deletion, and Section~\ref{sec:dualfix} describes how
\quads\ uses the excluded assignments.

\paragraph{Gilmore--Lawler bound.} The GLB~\cite{gilmore1962optimal,lawler1963quadratic}
bounds the quadratic cost of each allowed assignment from below by a linear assignment problem. For $M\in\S^N$ and $p=(i,j)\in\cE$, let
\begin{equation}\label{eq:H}
H_{ij}(M)\coloneqq M_{pp}+M_{0p}+M_{p0}+\min_{\sigma}\sum_{k\neq i}M_{p,(k,\sigma(k))},
\end{equation}
where $\sigma$ ranges over the bijections from $[m]\setminus\{i\}$ to $[m]\setminus\{j\}$. Thus $H_{ij}(M)$ is the
smallest contribution of the assignment $p$ to $y^\top My$ over the permutations with $\pi(i)=j$, and it requires a
linear assignment problem of order $m-1$. The GLB then solves a linear assignment problem of order $m$ with $H(M)$ as the cost matrix.
Let
\begin{equation}\label{eq:GL2}
\GL(M)\coloneqq M_{00}+\LAP(H(M)),
\end{equation}
with $\LAP$ defined in~\eqref{eq:lap}.
The GLB of the node is $c_0+\GL(L)$, and it costs $m^2$ small linear assignment problems.

\paragraph{Level-1 RLT bound.}\label{sec:rlt1} The level-1 RLT relaxation of~\eqref{eq:Y}~\cite{adams1994improved} is a
linear program, and its value is at least the GLB. We bound its value from below by the dual ascent of Hahn and
Grant~\cite{hahn1998lower}. Each round computes $c_0+\GL(M)$ for a lifted cost matrix $M$ with $y^\top My=y^\top Ly$ for the lift $y$ of every
permutation. The dual solutions of its linear assignment problems give the matrix $M$ of the next round. The $m^2$
linear assignment problems of a round are independent, and we solve them in parallel on the CPU cores. \quads\ computes
this bound only at the root and only if the GLB of the root is close to the objective value of the best known
permutation. Since the data are integers, we verify the final bound in exact integer arithmetic.

\paragraph{Quadratic programming bound.} The QPB~\cite{anstreicher2001new} is the minimum of a convex quadratic function over the
doubly stochastic matrices. This function equals the objective of~\eqref{eq:node} at every permutation matrix, so its minimum is a lower bound on $v^{\rm node}$.
\quads\ computes it by the Frank--Wolfe method~\cite{frank1956algorithm,brixius2001solving}. At every iterate, a feasible dual solution of the linearized problem gives a valid lower bound and reduced costs. The QPB requires symmetric $A$ and $B$. If only one of them is asymmetric, we replace it by its symmetric part, which does not change the objective.

\paragraph{DNN bound.} We relax~\eqref{eq:Y} by the SDP-RLT relaxation~\cite{hou2025rinnalplus}.
Its constraints are written in terms of the lifted matrix
\[
Y=\begin{pmatrix}Y_{00} & x_Y^\top\\ x_Y & {X}_Y\end{pmatrix}\in\S^N,\qquad x_Y\in\R^{m^2},\ {X}_Y\in\S^{m^2},
\]
which equals $yy^\top$ for a feasible $x$. The products of the constraints $\mathcal{A}x=e$ with $x$ give
$\mathcal{A}xx^\top=e\,x^\top$, and the binary constraints give $x_p^2=x_p$. The products of the valid inequalities
$e-x\ge0$ and $x\ge0$ with each other give $(e-x)x^\top\ge0$, $(e-x)(e-x)^\top\ge0$, and $xx^\top\ge0$. Relaxing $Y=yy^\top$ to $Y\succeq0$ gives
\begin{equation}\label{eq:sdprlt}
v^{\rm SDP\text{-}RLT}\coloneqq\min\left\{\inprod{L}{Y}+c_0\;\middle|\;\begin{array}{@{}l@{}}
Y_{00}=1,\ \ \mathcal{A}x_Y=e,\ \ \mathcal{A}{X}_Y=e\,x_Y^\top,\ \ \diag({X}_Y)=x_Y,\\
Y\in\S^N_+,\ \ Y\ge0,\\
e-x_Y\ge0,\ \ e\,x_Y^\top-{X}_Y\ge0,\ \ ee^\top-e\,x_Y^\top-x_Ye^\top+{X}_Y\ge0\end{array}\right\}.
\tag{SDP-RLT}
\end{equation}
For the QAP, the constraints in the third row are implied by the others~\cite[Proposition~8.2]{sherali1999reformulation},~\cite[Lemma~3]{hu2021linearization}, and so are the upper bounds $Y\le1$~\cite{oliveira2018admm}. Removing the third row of~\eqref{eq:sdprlt} therefore does not change its feasible set, so that
$v^{\rm SDP\text{-}RLT}=v^{\rm DNN}$ for
\begin{equation}\label{eq:dnn}
v^{\rm DNN}\coloneqq\min\Big\{\inprod{L}{Y}+c_0\;\Big|\;Y_{00}=1,\ \mathcal{A}x_Y=e,\ \mathcal{A}{X}_Y=e\,x_Y^\top,\
\diag({X}_Y)=x_Y,\ Y\in\S^N_+,\ Y\ge0\Big\}.\tag{DNN}
\end{equation}
Hence the level-1 RLT inequalities do not strengthen the DNN relaxation of the QAP. Since $Y$ is both positive semidefinite and nonnegative, we call~\eqref{eq:dnn} the DNN relaxation. It is equivalent to the DNN relaxations of the QAP in~\cite{zhao1998semidefinite,povh2009copositive}, and $v^{\rm DNN}\le v^{\rm node}$.

\paragraph{Choice of the bounds.} The QPB and the DNN bound cost much more than the GLB. We therefore compute them only at the nodes
that they are likely to prune, as estimated from their increases at the nodes of the same depth processed before. For every depth, \quads\ records the
increase of the lower bound achieved by the QPB and the increase of the DNN bound over the GLB at the nodes where these
bounds are computed, and the expected increase of a bound at a depth is the mean of the recorded increases. At a node,
a bound is computed only if a fixed multiple of its expected increase at the depth of the node reaches the gap between
the pruning threshold $\theta$ of Section~\ref{sec:alg} and the current lower bound of the node. Otherwise the node is branched without this bound. Both bounds
are computed at the first nodes of every depth, and the QPB also at a fixed fraction of the other nodes, so that the statistics remain
up to date.

\subsection{The algorithm}\label{sec:alg}

Algorithm~\ref{alg:bnb} summarizes the framework of \quads. Its lower bounds are those of Section~\ref{sec:dnn}. The DNN bound is computed by the
ADMM of Section~\ref{sec:bounds} with the stopping rule of Section~\ref{sec:stop}. Its cuts, reductions, dual
fixing, and branching rule are those of Section~\ref{sec:reduce}. The rest of this section describes the search, the pruning, and the upper bounds.

\paragraph{Search.}\label{sec:search} The best permutation found so far is called the incumbent, and $\bar v$ denotes its objective value. The nodes that have been created but not yet processed are called open. \quads\ processes them in increasing order of their lower bounds, and among equal bounds it prefers the node with more fixed assignments. The search ends when no node is open, and the incumbent is then optimal.

\paragraph{Pruning.} A node is pruned if it is discarded without branching. For integer data, all objective values lie in a set $v_0+\delta\mathbb{Z}$, where $v_0$ is the objective value of any permutation and $\delta\ge1$ is the greatest common divisor of the differences of the products $A_{ik}B_{jl}$. Hence a node can improve the incumbent only if its optimal value is at most the pruning threshold $\theta$, the largest element of $v_0+\delta\mathbb{Z}$ below $\bar v$. A node is pruned as soon as a valid lower bound exceeds $\theta$. A node is also pruned if it contains no permutation (Section~\ref{sec:screen}).

\paragraph{Upper bounds.}\label{sec:ub} Since the bounds prune a node only when its lower bound exceeds $\theta$, a good incumbent is needed early. We compute the initial upper bound at the root by robust tabu search~\cite{taillard1991robust}, a local search that swaps the locations of two facilities. During the tree search, every node gives two further candidate permutations. The first is the permutation $\pi$ of the node that maximizes $\sum_{i}(x_Y)_{(i,\pi(i))}$, where $x_Y$ is the block of the final $Y$ of the node. It is computed by a linear assignment problem. The second is the permutation that attains the GLB of the node. We improve both candidates by the swaps of the tabu search. Nodes with at most three free facilities are solved by enumeration. Every candidate is evaluated in integer arithmetic before it replaces the incumbent.

\begin{algorithm}[!htbp]
\caption{Branch and bound in \quads\ for solving~\eqref{eq:qap}}\label{alg:bnb}
\begin{algorithmic}[1]
\Require data $(A,B)$, optionally an upper bound supplied by the user
\If{only one of $A$ and $B$ is asymmetric} replace it by its symmetric part \Comment{Sec.~\ref{sec:dnn}} \EndIf
\State compute the initial incumbent $\bar{v}$  by tabu search \Comment{Sec.~\ref{sec:ub}}
\State compute the automorphisms of the data \Comment{Sec.~\ref{sec:exact}}
\State stop if the clique bound equals $\bar v$ \Comment{Sec.~\ref{sec:clique}}
\State put the root into the pool $\cP$
\While{$\cP\neq\emptyset$}
  \algphase{Select a node}
  \State remove a node with the smallest lower bound from $\cP$, and update $\theta$
  \If{the node is a duplicate or dominated} prune it and continue \Comment{Sec.~\ref{sec:exact}} \EndIf
  \If{at most three facilities are free or $A$ or $B$ is diagonal}
    \State solve the node exactly, update $\bar v$, and continue \Comment{Secs.~\ref{sec:ub}, \ref{sec:exact}}
  \EndIf
  \algphase{GLB, QPB, and level-1 RLT bound}
  \State $h\gets$ the larger of the GLB and, where it is computed, the QPB \Comment{Sec.~\ref{sec:dnn}}
  \If{$h>\theta$} prune the node and continue \EndIf
  \If{the node is the root and the GLB is close to $\bar v$}
    \State $h\gets\max\{h,\text{level-1 RLT bound}\}$, and prune the node if $h>\theta$ \Comment{Sec.~\ref{sec:rlt1}}
  \EndIf
  \algphase{DNN bound}
  \If{the DNN bound is likely to prune the node} \Comment{Sec.~\ref{sec:dnn}}
  \Repeat
    \State perform ADMM iterations, in block form if it is cheaper \Comment{Secs.~\ref{sec:admm}, \ref{sec:warm}, \ref{sec:clique}}
    \State $h\gets\max\{h,\text{lower bound~\eqref{eq:phi}}\}$ \Comment{Sec.~\ref{sec:cert}}
  \Until{$h>\theta$ or the stopping rule holds} \Comment{Sec.~\ref{sec:stop}}
  \EndIf
  \If{the node is the root} add aggregated RLT cuts and solve again \Comment{Sec.~\ref{sec:cuts}} \EndIf
  \algphase{Incumbent and dual fixing}
  \State update $\bar v$ by rounding the final $Y$ and a local search \Comment{Sec.~\ref{sec:ub}}
  \State exclude and fix assignments by dual fixing, and raise $h$ \Comment{Sec.~\ref{sec:dualfix}}
  \algphase{Prune or branch}
  \If{$h>\theta$ or the node contains no permutation} prune the node
  \Else
    \State branch on a line with one child per orbit, or on one assignment \Comment{Secs.~\ref{sec:exact}, \ref{sec:branch}}
    \State put the children with their warm starts into $\cP$
  \EndIf
\EndWhile
\Ensure an optimal permutation and the optimal value $v^{\rm QAP}$
\end{algorithmic}
\end{algorithm}

\section{The node solver}\label{sec:bounds}

We show how \quads\ solves the DNN relaxation of Section~\ref{sec:dnn} at a low cost. We also show how it
obtains a valid lower bound from an inexact solution. These two steps address the first two difficulties of
Section~\ref{sec:challenges}.

\subsection{Two-block ADMM}\label{sec:admm}

We separate the semidefinite constraint together with the
assignment equalities from the polyhedral constraints~\cite{hou2025rinnal}. The constraints of~\eqref{eq:dnn} imply $({X}_Y)_{pq}=0$ for every pair of assignments in
\[
\cG\coloneqq\big\{((i,j),(k,l))\in\cE\times\cE:\ i=k,\ j\neq l\ \text{ or }\ i\neq k,\ j=l\big\}.
\]
Indeed, let $p=(i,j)$. The equality of $\mathcal{A}{X}_Y=e\,x_Y^\top$ for the row sum of facility $i$ in the
column $p$ is $\sum_l({X}_Y)_{(i,l),p}=(x_Y)_p$. Since $({X}_Y)_{pp}=(x_Y)_p$, the terms with $l\neq j$ sum to zero,
and they are nonnegative. Hence $({X}_Y)_{(i,l),p}=0$ for $l\neq j$. The pairs with $i\neq k$ and $j=l$ follow in the same
way from the column sums. For the lift $Y=yy^\top$, these entries are $X_{ij}X_{kl}=0$.
These gangster
constraints~\cite{zhao1998semidefinite} state that two assignments that share a facility or a location cannot both hold. We keep them in the polyhedral part, where their projection costs nothing and sets these entries of every iterate to zero. Such redundant constraints in the polyhedral subproblem are also exploited in~\cite{oliveira2018admm}.
This gives the two closed convex sets
\begin{align*}
\cK_1&\coloneqq\{W\in\S^N_+:\ \mathcal{A}x_W=W_{00}e,\ \mathcal{A}{X}_W=e\,x_W^\top\},\\
\cK_2&\coloneqq\{Y\in\S^N:\ Y\ge0,\ Y_{00}=1,\ ({X}_Y)_\cG=0,\ Y_{pp}=Y_{0p}\ \ \forall\,p\in\cE\},
\end{align*}
where $W$ is partitioned as $Y$. Thus $\cK_1$ contains the lifted assignment equalities of~\eqref{eq:dnn} in
homogeneous form. The set $\cK_2$ contains the normalization, the
gangster constraints, and the arrow equalities $Y_{pp}=Y_{0p}$, which lift the binary constraints $x_p^2=x_p$. Hence
$\cK_1\cap\cK_2$ is the feasible set
of~\eqref{eq:dnn} with the implied gangster constraints, and~\eqref{eq:dnn} is equivalent to
\begin{equation}\label{eq:split}
\min_{Y,W\in\mathbb{S}^N} \Big\{\inprod{L}{Y}+c_0\;\Big|\;Y\in\cK_2,\ W\in\cK_1,\ Y=W\Big\}.
\end{equation}
We solve~\eqref{eq:split} by ADMM with multiplier $Z\in\S^N$ for the constraint $Y=W$, penalty parameter $\beta>0$,
and dual step length $\gamma=1.618$, which lies in the interval $(0,(1+\sqrt5)/2)$ for the convergence of ADMM established in~\cite{glowinski1975approximation,fazel2013hankel,yang2026inexact}. Starting from $(Y^0,Z^0)$, the $k$th iteration is given by
\begin{equation}\label{eq:admm}
\begin{aligned}
W^{k+1}&=\Proj_{\cK_1}\big(Y^k+Z^k/\beta\big),\\
Y^{k+1}&=\Proj_{\cK_2}\big(W^{k+1}-(L+Z^k)/\beta\big),\\
Z^{k+1}&=Z^k+\gamma\beta\,(Y^{k+1}-W^{k+1}).
\end{aligned}
\end{equation}
Let $P$ be the orthogonal projector onto the subspace $\cF\coloneqq\{(\tau;x)\in\R^N:\ \mathcal{A}x=\tau e\}$, which
contains the lift of every permutation. With $J\coloneqq I_m-ee^\top/m$ and $\hat e\coloneqq(1;e/m)\in\R^N$, it is given by
\[
P=\tfrac12\,\hat e\hat e^\top+\begin{pmatrix}0&0\\0&J\otimes J\end{pmatrix}.
\]
The dual problem of~\eqref{eq:split} is
\begin{equation}\label{eq:dual}
\max_{Z\in\mathbb{S}^N} \left\{c_0+(L+Z)_{00}\;\middle|\;\begin{array}{@{}l@{}}PZP\preceq0,\ \ (L+Z)_{pp}+2(L+Z)_{0p}\ge0\ \ \forall\,p\in\cE,\\[5pt]
(L+Z)_{pq}\ge0\ \ \forall\,(p,q)\in(\cE\times\cE)\setminus\cG,\ p\neq q\end{array}\right\}.
\end{equation}
The assignment equalities enter ADMM only through $P$, and both projections in~\eqref{eq:admm} have closed
forms. For $M\in\S^N$, let $M_+$ and $M_-$ be the parts of an eigenvalue decomposition of $M$ with the positive and the
negative eigenvalues, so that $M=M_++M_-$.

\begin{proposition}[Efficient projections in ADMM]\label{prop:proj}
Let $M\in\S^N$.
\begin{enumerate}[label=(\alph*),leftmargin=*,itemsep=1pt]
\item The projection onto $\cK_1$ is
\[
\Proj_{\cK_1}(M)=(PMP)_+.
\]
\item The projection onto $\cK_2$ is given entrywise by
\[
\big(\Proj_{\cK_2}(M)\big)_{ab}=\begin{cases}
1, & (a,b)=(0,0),\\
\max\{0,(M_{pp}+M_{0p}+M_{p0})/3\}, & (a,b)\in\{(p,p),(0,p),(p,0)\},\ p\in\cE,\\
0, & (a,b)\in\cG,\\
\max\{0,M_{ab}\}, & \textup{otherwise}.
\end{cases}
\]
\end{enumerate}
Moreover, $Pv$ costs $\cO(N)$ operations for $v\in\R^N$, and $PMP$ costs $\cO(N^2)$ operations.
\end{proposition}
\begin{proof}
(a) The equalities in $\cK_1$ state that every column of $W$ lies in $\cF$, so that
$\cK_1=\{W\in\S^N_+\mid (I-P)W=0\}=\{W\in\S^N_+\mid \inprod{I-P}{W}=0\}$, where the second equality holds since $I-P\succeq 0$ and $W\succeq0$. Hence (a) follows from~\cite[Lemma~2]{hou2025rinnalplus}. The formula of $P$ follows since $\cF$ is the orthogonal sum of the line spanned by $\hat e$ and of the vectors $(0;x)$ with $Xe=X^\top e=0$, onto which the projector is $J\otimes J$. Since
$JXJ=X-Xee^\top/m-ee^\top X/m+(e^\top Xe/m^2)ee^\top$, $Pv$ costs $\cO(m^2)=\cO(N)$
operations, and $PMP$ applies $P$ to the $N$ columns and rows of $M$.

(b) The set $\cK_2$ is a product of constraints on single entries and on the triples $(Y_{pp},Y_{0p},Y_{p0})$, which are
equal and nonnegative. Hence $\Proj_{\cK_2}$ clips single entries and replaces each triple by its clipped mean.
\end{proof}

\begin{remark}[Partial eigenvalue decomposition]\label{rem:eig} A full eigenvalue decomposition of $S\coloneqq PMP$ costs $\cO(N^3)$ operations and dominates the cost of an
iteration. However, only a few eigenvalues of $S$ change sign from one iteration to the next. Hence \quads\ computes
only the eigenpairs of the side of the spectrum with fewer eigenvalues, using $S_+=S-S_-$ when the negative side has fewer eigenvalues, and it tracks this eigenspace across iterations by a Chebyshev-filtered subspace
iteration~\cite{li2020lowrank}. A full decomposition is computed only periodically or when the tracking fails.
\end{remark}

\paragraph{Penalty parameter.}
The default penalty parameter of a node is $\beta_0=\omega\max(1,\|L\|/N)$ with a constant $\omega>0$. During ADMM, we compare the primal residual $\eta_p=\|Y^k-W^k\|/(1+\|Y^k\|)$ with the dual residual $\eta_d=\beta\|Y^k-Y^{k-1}\|/(1+\|L\|)$. When their ratio stays away from a target ratio over several checks, $\beta$ is multiplied or divided by a constant factor, within a fixed interval around $\beta_0$~\cite{yang2015sdpnalplus}. The initial value of $\beta$ at the children is given in Section~\ref{sec:warm}.

\subsection{Valid lower bounds}\label{sec:cert}

ADMM is stopped long before convergence, so its multiplier $Z$ is in general infeasible for the dual~\eqref{eq:dual},
and the dual objective value $c_0+(L+Z)_{00}$ is not a valid lower bound on $v^{\rm node}$. We obtain a valid lower bound from every
$Z\in\S^N$ as the dual function of another relaxation of~\eqref{eq:node}. Let
\begin{equation}\label{eq:Omega}
\Omega\coloneqq\left\{Y=\begin{pmatrix}1&x_Y^\top\\ x_Y&X_Y\end{pmatrix}\in\R^{N\times N}\;\middle|\;
\begin{array}{@{}l@{}}\mathcal{A}x_Y=e,\ \ X_Y\mathcal{A}^\top=x_Ye^\top,\ \ x_Y\ge0,\ \ X_Y\ge0,\\[3pt] (X_Y)_{pq}=0\ \ \forall\,(p,q)\in\cG\end{array}\right\}.
\end{equation}
Every $y\in\cY$ in~\eqref{eq:Y} satisfies $yy^\top\in\Omega\cap\cK_1$ and $\|y\|^2=m+1$. Hence the problem
\begin{equation}\label{eq:tight}
\min\Big\{\inprod{L}{Y}+c_0\;\Big|\;Y\in\Omega,\ W\in\cK_1,\ \tr W=m+1,\ \sym(Y)=W\Big\}
\end{equation}
has an optimal value at most $v^{\rm node}$. Since $\inprod{Z}{\sym(Y)}=\inprod{Z}{Y}$ for $Z\in\S^N$, its Lagrangian
dual function is
\begin{equation}\label{eq:dfun}
d(Z)\coloneqq c_0+\min_{Y\in\Omega}\inprod{L+Z}{Y}+\min\big\{\inprod{-Z}{W}\;\big|\;W\in\cK_1,\ \tr W=m+1\big\},
\end{equation}
and weak duality gives $v^{\rm node}\ge d(Z)$ for every $Z\in\S^N$.

\begin{lemma}[Gilmore--Lawler bound as a linear program~\cite{adams1994improved}]\label{lem:GL}
For every $M\in\S^N$, $\min\{\inprod{M}{Y}\mid Y\in\Omega\}=\GL(M)$. Moreover, for every permutation $\pi$ with
lift $y$,
\begin{equation}\label{eq:GLbound}
y^\top My\ \ge\ M_{00}+\sum_{i=1}^mH_{i,\pi(i)}(M)\ \ge\ \GL(M).
\end{equation}
\end{lemma}
The set $\Omega$ is the feasible set of the level-1 RLT relaxation of Adams and Johnson~\cite{adams1994improved} without its symmetry
constraints $Y_{pq}=Y_{qp}$. We omit these constraints so that, by Lemma~\ref{lem:GL}, the first minimum in~\eqref{eq:dfun} is $\GL(L+Z)$, which costs $m^2$ small linear assignment problems. With those constraints, it would be the optimal value of a linear program with $\cO(m^4)$ variables, which is expensive to compute. 
Lagrangian relaxations of the symmetry constraints also give the dual-ascent bounds of~\cite{hahn1998lower,adams2007level2}. The second minimum in~\eqref{eq:dfun} is bounded by the largest eigenvalue of $PZP$, which gives the following bound.

\begin{proposition}[Lower bound from any multiplier]\label{prop:cert}
For every $Z\in\S^N$ and every $t\ge\lambda_{\max}(PZP)$, let $M\coloneqq L+Z-tI$. Then
\begin{equation}\label{eq:phi}
v^{\rm node}\ \ge\ d(Z)\ \ge\ \varphi(Z,t)\coloneqq c_0+\GL(M)=c_0+\GL(L+Z)-(m+1)\,t,
\end{equation}
and $d(Z)=\varphi(Z,t)$ if $t=\lambda_{\max}(PZP)>0$.
\end{proposition}
\begin{proof}
By Lemma~\ref{lem:GL}, the first minimum in~\eqref{eq:dfun} equals $\GL(L+Z)$. Every $W\in\cK_1$ satisfies $W=PWP$ and $W\succeq0$,
so that $\inprod{Z}{W}=\inprod{PZP}{W}\le(m+1)\lambda_{\max}(PZP)$ whenever $\tr W=m+1$. If $\lambda_{\max}(PZP)>0$,
equality holds for $W=(m+1)vv^\top$ with a unit eigenvector $v$ of $PZP$ for this eigenvalue, since then $v=Pv$ and
$W\in\cK_1$. Hence $d(Z)\ge c_0+\GL(L+Z)-(m+1)t$, with equality if $t=\lambda_{\max}(PZP)>0$. In particular, $W=yy^\top$ gives
$y^\top Zy\le(m+1)t$, so that $y^\top Ly\ge y^\top My$ for every $y\in\cY$. Subtracting $tI$ decreases
$M_{00}$ and every $M_{pp}$ by $t$ and leaves the other entries unchanged, so that $H(M)=H(L+Z)-t\,ee^\top$ and
$\GL(M)=\GL(L+Z)-(m+1)t$.
\end{proof}

Every lower bound that \quads\ computes from a multiplier is a lower estimate of the dual function
of~\eqref{eq:tight}, possibly with additional constraints. The bounds of dual fixing (Section~\ref{sec:dualfix}) and of the cuts (Section~\ref{sec:cuts}) are of this form. The
next theorem shows that the bound~\eqref{eq:phi} loses nothing at an optimal multiplier.

\begin{theorem}[Tightness of the bound]\label{thm:ceiling}
Let $\varphi_0(Z)\coloneqq\varphi(Z,\lambda_{\max}(PZP))$. Then the following hold.
\begin{enumerate}[label=(\alph*),leftmargin=*,itemsep=1pt]
\item $\varphi_0(Z)\le v^{\rm DNN}$ for every $Z\in\S^N$.
\item $\varphi_0(Z_*)=v^{\rm DNN}$ for every optimal solution $Z_*$ of~\eqref{eq:dual}.
\end{enumerate}
\end{theorem}
\begin{proof}
(a) Let $Y$ be feasible for~\eqref{eq:dnn}. Then $Y\in\Omega$, since $X_Y\mathcal{A}^\top=x_Ye^\top$ is the
transpose of $\mathcal{A}X_Y=e\,x_Y^\top$ and the gangster constraints hold. Moreover, $Y\in\cK_1$ since its columns lie
in $\cF$, and $\tr Y=Y_{00}+e^\top x_Y=m+1$. Hence $(Y,Y)$ is feasible for~\eqref{eq:tight}, and Proposition~\ref{prop:cert} and weak duality give $\varphi_0(Z)\le d(Z)\le c_0+\inprod{L}{Y}$. Minimizing over $Y$ gives
$\varphi_0(Z)\le v^{\rm DNN}$.

(b) Let $Z_*$ be an optimal
solution of~\eqref{eq:dual}, which exists and has the optimal value $v^{\rm DNN}$ by Fenchel's duality theorem~\cite[Theorem~31.1]{rockafellar1970convex} since $\cK_2$ is polyhedral and the barycenter of the matrices $yy^\top$ over the lifts $y$ of all
permutations lies in the relative interior of $\cK_1$. Dual feasibility gives $PZ_* P\preceq0$. Since $PZ_*P$ is singular, $\lambda_{\max}(PZ_*P)=0$ and hence $t=0$. Moreover,
$H(L+Z_*)\ge0$ by the constraints of~\eqref{eq:dual}, since the entries $(p,(k,\sigma(k)))$ with $p=(i,j)$ in~\eqref{eq:H} do not
lie in $\cG$. Hence
$\varphi_0(Z_*)\ge c_0+(L+Z_*)_{00}=v^{\rm DNN}$, and equality holds by (a).
\end{proof}

Every fixed number of ADMM iterations, \quads\ evaluates~\eqref{eq:phi} at the multiplier
$Z=Z^k+\beta(Y^k-W^{k+1})$. It keeps the largest value $h$ as the lower bound of the node, with a margin for rounding errors (Section~\ref{sec:verify}). By Theorem~\ref{thm:ceiling},
a node with $v^{\rm DNN}\le\theta$ cannot be pruned in this way, however long ADMM runs. In such a case, the node
should therefore be branched early, and the stopping rule of Section~\ref{sec:stop} detects it.

\subsection{Early termination}\label{sec:stop}

A node is pruned at the first evaluation at which its best
lower bound $h_k$ after $k$ ADMM iterations exceeds
$\theta$. If $h_k$ stalls below $\theta$, further iterations are wasted and the node should be branched. We
therefore stop the computation of the bound when a forecast of its growth does
not reach~$\theta$~\cite{rendl2010solving,krislock2014improved}. We forecast the limit of the bound by a
geometric extrapolation. Let $w$ be a fixed number of iterations. At every evaluation with $k\ge2w$, let
$\Delta_1=\max(0,h_k-h_{k-w})$, $\Delta_2=\max(0,h_{k-w}-h_{k-2w})$, and
\[
\kappa=\begin{cases}\min\{\Delta_1/\Delta_2,\bar\kappa\} & \text{if } \Delta_2>0,\\ \bar\kappa & \text{if } \Delta_2=0<\Delta_1,\\ 0 & \text{if } \Delta_1=\Delta_2=0,\end{cases}
\]
with a cap $\bar\kappa<1$. If the increments of the bound decay geometrically with ratio $\kappa$, then the bound converges to
$h_k+\Delta_1\kappa/(1-\kappa)$, which, for $\kappa=\Delta_1/\Delta_2$, is Aitken's $\Delta^2$ extrapolation of
$h_{k-2w},h_{k-w},h_k$. The forecast fails when $h_k+\xi\,\Delta_1\kappa/(1-\kappa)\le\theta$, where the factor $\xi>1$
makes the forecast optimistic. If $\theta-h_k$ exceeds $\tau_{\rm gap}\max(1,|\theta|)$ for a constant $\tau_{\rm gap}>0$, the node stops at the first failure. Otherwise it stops at the second failure. At the root, the rule requires several consecutive failures, and it is not used for large $m$. A node also stops after a maximal number of iterations or a maximal time. A
stopped node is branched. Hence the rule affects the running time but not the correctness of the result.

\subsection{Warm start of children}\label{sec:warm}
ADMM can start from any pair $(Y,Z)$. Each child therefore starts from the final variables of its parent.
Let a child fix the assignment $a=(i,j)$, and let the subscript ${\rm c}$ mark the quantities of
Section~\ref{sec:admm} for the child, whose lifted matrices have order $N_{\rm c}\coloneqq(m-1)^2+1$. Every
$y\in\cY$ in~\eqref{eq:Y} with $y_a=1$ satisfies $y=Ty_{\rm c}$ for the lift $y_{\rm c}\in\cY_{\rm c}$ of the same permutation in the child. Here $T\in\{0,1\}^{N\times N_{\rm c}}$ copies the entry $0$ of $y_{\rm c}$ to the entries $0$
and $a$ of $y$ and every other entry to the entry of the same assignment. The free facilities and locations of the child keep their indices in the parent. Hence $y_{\rm c}^\top T^\top LTy_{\rm c}=y^\top Ly$. Thus $T^\top LT$ represents the objective of the child
and differs from its matrix~\eqref{eq:L} only by a constant on $\cK_2^{\rm c}$, by the arrow equalities and $(Y_{\rm c})_{00}=1$. The child starts from
\begin{equation}\label{eq:warm}
Z^0_{\rm c}=T^\top\hat ZT,\qquad Y^0_{\rm c}=\Proj_{\cK_2^{\rm c}}(\hat Y_{\rm c}),
\end{equation}
where $\hat Z$ is the multiplier of the parent with the largest value of~\eqref{eq:phi} and $\hat Y_{\rm c}$ is the principal
submatrix of the final $Y$ of the parent on the indices of the child. The initial penalty parameter of the child is its default value multiplied by the ratio of the final to the
default penalty parameter of the parent. The next proposition justifies~\eqref{eq:warm}.

\begin{proposition}[Exact warm start]\label{prop:warm}
Let $(Y^*,Z^*)$ be an optimal primal--dual pair of~\eqref{eq:split} and~\eqref{eq:dual} at the parent with
$Y^*_{aa}=1$, and let $\hat Y_{\rm c}$ be the principal submatrix of $Y^*$ on the indices of the child. Then the
following hold.
\begin{enumerate}[label=(\alph*),leftmargin=*,itemsep=1pt]
\item $Y^*=T\hat Y_{\rm c}T^\top$.
\item $(\hat Y_{\rm c},T^\top Z^*T)$ is an optimal primal--dual pair of the child, and the DNN bounds of the child and of
the parent coincide.
\item This pair is a fixed point of~\eqref{eq:admm} for the child, and it is the start~\eqref{eq:warm} with $\hat Z=Z^*$.
\end{enumerate}
\end{proposition}
\begin{proof}
(a) Since $Y^*\succeq0$ and $Y^*_{00}=Y^*_{0a}=Y^*_{aa}=1$, the rows $0$ and $a$ of $Y^*$ are equal. The other
assignments with the facility $i$ or the location $j$ have $Y^*_{pp}=0$, so their rows are zero. Hence
$Y^*=T\hat Y_{\rm c}T^\top$.

(b) The pair $(Y,Z)$ is optimal for~\eqref{eq:split} and~\eqref{eq:dual} if and only if it satisfies the
Karush--Kuhn--Tucker (KKT) conditions
\[
Y\in\cK_1\cap\cK_2,\qquad Z\in N_{\cK_1}(Y),\qquad -(L+Z)\in N_{\cK_2}(Y),
\]
where $N_\cK(Y)$ is the normal cone of $\cK$ at $Y$. By (a), $\hat Y_{\rm c}\in\cK_1^{\rm c}\cap\cK_2^{\rm c}$. Moreover,
$TY_{\rm c}T^\top\in\cK_l$ for every $Y_{\rm c}\in\cK_l^{\rm c}$, $l=1,2$, so that
\[
\inprod{T^\top Z^*T}{Y_{\rm c}-\hat Y_{\rm c}}=\inprod{Z^*}{TY_{\rm c}T^\top-Y^*}\le0\qquad\forall\,Y_{\rm c}\in\cK_1^{\rm c},
\]
and hence $T^\top Z^*T\in N_{\cK_1^{\rm c}}(\hat Y_{\rm c})$. In the same way, $-T^\top(L+Z^*)T\in N_{\cK_2^{\rm c}}(\hat
Y_{\rm c})$. Since $T^\top LT$ represents the objective of the child, the KKT conditions of the child hold at $(\hat Y_{\rm
c},T^\top Z^*T)$, and $\inprod{T^\top LT}{\hat Y_{\rm c}}=\inprod{L}{Y^*}$.

(c) The KKT conditions state that the pair is a fixed point of~\eqref{eq:admm}. Since $\hat Y_{\rm c}\in\cK_2^{\rm c}$,
$\Proj_{\cK_2^{\rm c}}(\hat Y_{\rm c})=\hat Y_{\rm c}$, so that the pair is the start~\eqref{eq:warm} with $\hat Z=Z^*$.
\end{proof}

\section{Reducing the search tree}\label{sec:reduce}

We reduce the number of nodes of the search tree by dual fixing, cuts at the root, the structure of the data,
and the branching rule. These address the third difficulty of Section~\ref{sec:challenges}. The
structure of the data also reduces the cost of a node by block diagonalization and, in the extreme case, replaces the
tree by an enumeration.

\subsection{Dual fixing}\label{sec:dualfix}

Dual fixing can discard children of a node without solving them. From the dual information of the parent node, in
particular the multiplier of its DNN bound, \quads\ bounds from below the optimal values of the two children with $X_a=1$ and
$X_a=0$ of every assignment $a$. If the bound of one child exceeds the pruning threshold $\theta$ of Section~\ref{sec:search}, then $X_a$ is fixed to the value of the other
child in the whole subtree of the parent node. This is the
analogue of reduced-cost fixing in mixed-integer programming, in which the reduced costs of the relaxation of a node
fix variables without solving the corresponding subproblems. For an allowed assignment $a$, let $v^1_a$ and
$v^0_a$ be the optimal values of the children with $X_a=1$ and $X_a=0$, which are the smallest objective values of
the permutations of the node that use $a$ and that avoid $a$, respectively. Given lower bounds $\underline v^1_a\le v^1_a$
and $\underline v^0_a\le v^0_a$, \quads\ applies the rule
\begin{equation}\label{eq:dfrule}
X_a=0 \ \text{ if } \underline v^1_a>\theta,\qquad X_a=1 \ \text{ if } \underline v^0_a>\theta,
\end{equation}
in the subtree of the parent node. The rule is valid since the discarded child contains no permutation with objective value
at most $\theta$. 
The first case is called the exclusion of $a$, while the second case is the fixing of $a$.

\paragraph{Bounds from the multiplier.} We propose to obtain $\underline v^1_a$ from the multiplier of the
DNN bound of the node. Restricting the minimum over the permutation matrices in $\varphi(Z,t)$ of
Proposition~\ref{prop:cert} to $X_{ij}=1$ or to $X_{ij}=0$ corresponds to adding
the constraint $(x_Y)_a=1$ or $(x_Y)_a=0$ to $\Omega$ in~\eqref{eq:tight} and gives the following bounds.

\begin{proposition}[Assignment bounds from a multiplier]\label{prop:fix}
Let $Z\in\S^N$, $t\ge\lambda_{\max}(PZP)$, and $M\coloneqq L+Z-tI$. For every assignment $a=(i,j)\in\cE$, define
\begin{align}
f^1_a(Z,t)&\coloneqq c_0+M_{00}+\min\{\inprod{H(M)}{X}\mid X\in\Pi_m,\ X_{ij}=1\},\label{eq:fix}\\
f^0_a(Z,t)&\coloneqq c_0+M_{00}+\min\{\inprod{H(M)}{X}\mid X\in\Pi_m,\ X_{ij}=0\}.\label{eq:fixbar}
\end{align}
Then the matrix $X$ of every permutation of the node with objective value at most $\theta$ satisfies
$X_{ij}=0$ if $f^1_a(Z,t)>\theta$ and $X_{ij}=1$ if $f^0_a(Z,t)>\theta$.
\end{proposition}
\begin{proof}
By the proof of Proposition~\ref{prop:cert}, $y^\top Ly\ge y^\top My$ for every $y\in\cY$. Hence,
by~\eqref{eq:GLbound}, every permutation of the node with permutation matrix $X$ has objective value at least
$c_0+M_{00}+\inprod{H(M)}{X}$. This value is at least $f^1_a(Z,t)$ if $X_{ij}=1$ and at least $f^0_a(Z,t)$ if
$X_{ij}=0$, which gives the two rules.
\end{proof}

\begin{corollary}\label{cor:fix}
For every assignment $(i,j)\in\cE$, $f^0_{ij}(Z,t)=\min_{l\neq j}f^1_{il}(Z,t)=\min_{k\neq i}f^1_{kj}(Z,t)$. Hence, if
$f^0_{ij}(Z,t)>\theta$, then the first rule of Proposition~\ref{prop:fix} excludes every other assignment of the
facility $i$ and of the location $j$.
\end{corollary}
\begin{proof}
A permutation matrix with $X_{ij}=0$ has $X_{il}=1$ for exactly one $l\neq j$ and $X_{kj}=1$ for exactly one $k\neq i$.
\end{proof}

\quads\ computes $f^1_a\coloneqq f^1_a(\hat Z,t)$ at the multiplier $\hat Z$ of the best lower bound of the node
and its value of $t$. In addition, the child with $X_a=1$ receives the transferred multiplier $T^\top\hat ZT$ of
Section~\ref{sec:warm}, and its GLB at this multiplier,
\[
g^1_a\coloneqq c_0+\GL\big(T^\top(L+\hat Z)T\big)-(m+1)\,t,
\]
bounds the permutations that use $a$, i.e., $X_a=1$. For the lift $y_{\rm c}$ of such a permutation in the child, $y=Ty_{\rm c}$ is a lift of the node, so that $Py=y$, $\|y\|^2=m+1$, and $y_{\rm c}^\top T^\top\hat ZTy_{\rm c}=y^\top\hat Zy\le(m+1)t$. Lemma~\ref{lem:GL} applied to $T^\top(L+\hat Z)T$ in the child gives the claim. The bounds $f^1_a$ and $g^1_a$ cost $m^2+2|\cD|$ linear assignment problems of order $m-1$ per node, of which $m^2$ give $H(M)$, and $m^2$ shortest-path computations on graphs with $m-1$ nodes, which reuse these problems for $g^1_a$. In total, they cost $\cO(m^5)$ operations.

\begin{remark}[Relation to other fixing rules]\label{rem:fixrel}
The bound $f^1_a(0,0)$ is a lower bound on the GLB of the node with $X_a=1$. The fixing by the
GLB and by the reduced costs of the QPB~\cite{anstreicher2001new,brixius2001solving} does not use
a multiplier of the DNN relaxation. Variable fixing from semidefinite relaxations has
been used before~\cite{helmberg2000fixing,hou2026scsdp}. To our knowledge, it has not
been used with the DNN relaxation in published branch-and-bound methods for the QAP.
\end{remark}

\paragraph{Exclusion and fixing.}\label{sec:screen} For the QPB of Section~\ref{sec:dnn}, the reduced costs of a Frank--Wolfe iterate give a lower bound $b^{\rm QPB}_a$ on the child with $X_a=1$ before this child is solved~\cite{brixius2001solving}. \quads\ uses the best of the bounds available at the node,
\[
\underline v^1_a\coloneqq\max\big\{f^1_a(0,0),\ b^{\rm QPB}_a,\ f^1_a,\ g^1_a,\ \underline v_a^{1,{\rm par}}\big\},
\]
where $f^1_a(0,0)$ is~\eqref{eq:fix} at $Z=0$ and $t=0$, $b^{\rm QPB}_a$ is the bound from the reduced costs of the QPB, and $\underline v_a^{1,{\rm par}}$ is the bound of the parent. A bound that is not computed at
the node is omitted. Every assignment with $\underline v^1_a>\theta$ is excluded. 
The exclusions are then propagated. Every permutation of the node assigns
each facility to a different allowed location, which is called a perfect matching of the allowed assignments. An
assignment that belongs to no such matching is excluded. A facility with a single allowed location is fixed to it. If no
perfect matching exists, the node contains no permutation and is pruned. Otherwise its lower bound is raised to $\max_i\min_{j}
\{ \underline
v^1_{ij} \mid (i,j)\in\cD\}$, since the permutations that use an excluded assignment have objective values above $\theta$. The children inherit the remaining allowed assignments and their bounds, which the branching rule of
Section~\ref{sec:branch} uses.

\paragraph{Smaller relaxations after exclusions.}\label{par:smallrel} The exclusions also make the DNN relaxation of the node smaller, since every excluded
assignment $p\not\in \mathcal{D}$ sets the entry $y_p=0$. They accumulate with the depth, and on some instances most assignments of a
deep node are excluded.
\quads\ deletes the rows and columns corresponding to the 
excluded assignments from the matrix variable in the DNN relaxation and hence solves a smaller problem.

\subsection{Aggregated RLT cuts}\label{sec:cuts}

On some instances, the DNN bound of the root is far below the optimal value,
and the tree becomes large. The SDP-RLT relaxation~\eqref{eq:sdprlt} contains the products of two constraints of the
assignment system, which form the level-1 RLT, and these products are implied by~\eqref{eq:dnn} (Section~\ref{sec:dnn}). We therefore strengthen the bound by the level-2 RLT, which multiplies three constraints~\cite{adams2007level2}. However, its lifted variables have $\Theta(m^6)$ entries, which is too many for a node.
We therefore sum the level-2 RLT over groups of facilities. We add the resulting inequalities at the root and
pass them to the whole tree through the lifted cost matrix.

Let $\mathcal P_1,\dots,\mathcal P_r$ be a partition of the facilities into groups. For two distinct locations
$j,j'$ and two sets $U,U'$ of facilities, let
\[
Y(U,j;\,U',j')\coloneqq\sum_{i\in U}\sum_{i'\in U'}Y_{(i,j),(i',j')}.
\]
For the lift of a permutation, $Y(U,j;\,U',j')$ is one if $j$ receives a facility of $U$ and $j'$ a facility of
$U'$, and zero otherwise. For three distinct locations $j_1,j_2,j_3$, let $z_{i_1i_2i_3}\ge0$ be the variable of the
level-2 RLT~\cite{adams2007level2} for the product $x_{i_1j_1}x_{i_2j_2}x_{i_3j_3}$. The products of the assignment
equalities of the locations with the products of two variables give
\begin{equation}\label{eq:rlt2eq}
Y_{(i_1,j_1),(i_2,j_2)}=\sum_{i_3}z_{i_1i_2i_3},\qquad Y_{(i_1,j_1),(i_3,j_3)}=\sum_{i_2}z_{i_1i_2i_3},\qquad
Y_{(i_2,j_2),(i_3,j_3)}=\sum_{i_1}z_{i_1i_2i_3}.
\end{equation}

\begin{proposition}\label{prop:aggrlt}
Let $j_1,j_2,j_3$ be distinct locations, and let $U_1,U_2,U_3$ be unions of groups, chosen independently. Let
$\bar U_1$ be the set of the facilities outside $U_1$. Then
\begin{eqnarray}\label{eq:rlt2sum}
&& \sum_{i_1\notin U_1,\ i_2\in U_2,\ i_3\notin U_3}z_{i_1i_2i_3}+\sum_{i_1\in U_1,\ i_2\notin U_2,\ i_3\in U_3}z_{i_1i_2i_3}
\\
&=&
Y(U_2,j_2;\,\bar U_1,j_1)+Y(U_1,j_1;\,U_3,j_3)-Y(U_2,j_2;\,U_3,j_3).
\nonumber
\end{eqnarray}
Hence every $Y$ that is the projection of a feasible point of the level-2 RLT relaxation satisfies
\begin{equation}\label{eq:aggrlt}
Y(U_2,j_2;\,U_3,j_3)\le Y(U_2,j_2;\,\bar U_1,j_1)+Y(U_1,j_1;\,U_3,j_3).
\end{equation}
\end{proposition}
\begin{proof}
By~\eqref{eq:rlt2eq}, we have
$Y(U_2,j_2;\bar U_1,j_1) = \sum_{i_2\in U_2}\sum_{i_1\not\in U_1} (\sum_{i_3\in U_3} z_{i_1i_2i_3} +\sum_{i_3\not\in U_3} z_{i_1i_2i_3}).$
By expanding $Y(U_1,j_1;U_3,j_3)$ and $Y(U_2,j_2;U_3,j_3)$
similarly, we get the identity \eqref{eq:rlt2sum}.
Since $z\ge0$, the left-hand side of~\eqref{eq:rlt2sum} is
nonnegative, and \eqref{eq:aggrlt} follows.
\end{proof}

The left-hand side of~\eqref{eq:rlt2sum} is a sum of nonnegativity constraints of the level-2 RLT. We choose them
so that their cubic variables cancel, and only entries of $Y$ remain on the right-hand side. Hence
\eqref{eq:aggrlt} is a linear inequality $\inprod{V}{Y}\le0$ with $V\in\S^N$, and the $\Theta(m^6)$ variables of
the level-2 RLT do not enter a node. These inequalities are ST-inequalities of J\"unger and
Kaibel~\cite{junger2001box}, and the grouped cuts of~\cite{faye2005cutting} are among them.

\paragraph{Separation.}
\quads\ partitions the facilities into $r$ groups, where $r$ grows with $m$, by $k$-medoids clustering with $A_{ik}$ as the dissimilarity
between the facilities $i$ and $k$. Each group then has small internal flow. \quads\ uses the cuts only if the flows within these groups are much smaller than those
between the groups. It first runs a fixed number of ADMM
iterations without cuts. The separation computes $Y(\mathcal P,j;\,\mathcal P',j')$ for all groups $\mathcal P,\mathcal P'$ and all locations $j\neq j'$ at the current $Y$, so that every term of~\eqref{eq:aggrlt} is a sum of these values. It then enumerates all
triples of locations and all choices of unions. It adds the most violated inequalities $\inprod{V_c}{Y}\le0$, $c\in\mathcal C$, to the set $\cK_2$ of
Section~\ref{sec:admm}. We keep a multiplier $\mu_c\ge0$ for each cut $c$, and let $L_\mu\coloneqq L+\sum_{c\in\mathcal C}\mu_cV_c$. The step for $Y$ in~\eqref{eq:admm} becomes
\begin{equation}\label{eq:ycut}
Y^{k+1}=\Proj_{\cK_2}\big(W^{k+1}-(L_\mu+Z^k)/\beta\big).
\end{equation}
Let $\cK_2^{\mathcal C}\coloneqq\{Y\in\cK_2:\ \inprod{V_c}{Y}\le0\ \ \forall\,c\in\mathcal C\}$. The projection of $W^{k+1}-(L+Z^k)/\beta$ onto $\cK_2^{\mathcal C}$ equals the right-hand side of~\eqref{eq:ycut} for multipliers $\mu\ge0$ that satisfy $\inprod{V_c}{Y^{k+1}}\le0$ and $\mu_c\inprod{V_c}{Y^{k+1}}=0$ for every $c\in\mathcal C$. We approximate these multipliers by one pass of coordinate ascent on the dual of this projection~\cite{hildreth1957quadratic}, started from the multipliers of the previous iteration. A similar projection step is used for graph partition in~\cite{demeijer2023partitioning}. The pass replaces each $\mu_c$ in turn by the smallest $\mu_c\ge0$ with $\inprod{V_c}{\Proj_{\cK_2}(W^{k+1}-(L_\mu+Z^k)/\beta)}\le0$.
The separation runs in a number of rounds that grows with $m$, each followed by a fixed number of ADMM iterations. The cuts are discarded if they close only a small fraction of the gap between $\theta$ and the root bound without cuts.

\begin{proposition}\label{prop:cuts}
For every $\mu\ge0$, $y^\top Ly\ge y^\top L_\mu y$ for the lift $y$ of every
permutation. Moreover, Proposition~\ref{prop:cert} with $L_\mu$ in place of $L$ gives a valid lower bound on $v^{\rm node}$.
\end{proposition}
\begin{proof}
Every cut is valid for every permutation, so that $\mu_cy^\top V_cy\le0$ and adding the constraints
$\inprod{V_c}{Y}\le0$ to~\eqref{eq:tight} keeps $yy^\top$ feasible for every $y\in\cY$. With the multipliers $\mu_c$ of
these constraints, the Lagrangian dual function of the resulting problem is~\eqref{eq:dfun} with $L_\mu$ in place of
$L$, and the proof of Proposition~\ref{prop:cert} applies to it.
\end{proof}

After the separation, \quads\ freezes $\mu$ and solves the root once more with the lifted cost matrix $L_\mu$. It then uses $L_\mu$ in
the whole tree without further separation. As in Section~\ref{sec:warm}, a child receives the lifted cost matrix $T^\top
L_\mu T$, and Proposition~\ref{prop:cuts} remains valid at every node since $y=Ty_{\rm c}$. The term
$\sum_c\mu_cV_c$ is added to $L$ with downward rounding, so that the bound remains valid in floating-point
arithmetic (Section~\ref{sec:verify}).

\subsection{Zero flow and automorphisms}\label{sec:exact}

\begin{figure}[!htbp]
\centering
\colorlet{cA}{dtime}\colorlet{cB}{dvalid}
\begin{tikzpicture}[font=\small,
  fac/.style={circle,draw,thick,minimum size=6.5mm,inner sep=0pt},
  loc/.style={rectangle,draw,thick,minimum size=6mm,inner sep=0pt},
  note/.style={font=\footnotesize,align=center},
  sub/.style={font=\footnotesize,gray!60!black}]
\node at (1.5,2.3) {Zero flow};
\node[sub] at (1.5,1.85) {flow \textcolor{cA}{$A$}};
\node[fac] (f1) at (0.6,1.1) {1};
\node[fac] (f2) at (2.4,1.1) {2};
\draw[very thick,cA] (f1) -- (f2);
\node[fac] (f3) at (0.6,0) {3};
\node[fac] (f4) at (2.4,0) {4};
\node[note,anchor=north] at (1.5,-0.7) {3 and 4: zero flow\\ 1 and 2 fixed $\Rightarrow$ an LAP};
\node at (5.6,2.3) {Orbital branching};
\node[sub] at (5.6,1.85) {location reversal $g$};
\foreach \j/\x in {1/4.25,2/5.15,3/6.05,4/6.95} \node[loc] (l\j) at (\x,0.15) {\j};
\draw[very thick,cB] (l1)--(l2)--(l3)--(l4);
\draw[<->,dashed,gray] (l1.north) to[out=45,in=135,looseness=0.7] node[above,font=\scriptsize,inner sep=1pt] {$g$} (l4.north);
\draw[<->,dashed,gray] (l2.north) to[out=60,in=120,looseness=0.8] (l3.north);
\node[font=\footnotesize,cB] at (5.6,-0.4) {$B_{jl}=|j-l|$};
\node[note,anchor=north] at (5.6,-0.7) {children at locations $j$ 
\\  and $g(j)$ have equal values};
\node at (9.8,2.3) {Block diagonalization};
\node[sub] at (9.8,1.85) {exchange of facilities 3 and 4};
\draw[thick] (8.2,-0.35) rectangle (9.4,0.95);
\node[font=\scriptsize] at (8.8,0.3) {$1+m^2$};
\draw[->,thick] (9.5,0.3) -- (9.9,0.3);
\draw[thick,fill=dtime!12] (10.0,-0.15) rectangle (11.1,0.95);
\node[font=\scriptsize] at (10.55,0.4) {$1+Km$};
\draw[gray,fill=gray!20] (11.1,-0.4) rectangle (11.35,-0.15);
\draw[gray,fill=gray!20] (11.35,-0.65) rectangle (11.6,-0.4);
\node[note,anchor=north] at (9.8,-0.7) {$K$ classes\\ $\Rightarrow$ order $1+Km$};
\node at (13.6,2.3) {Clique bound};
\node[sub] at (13.6,1.85) {$Q=\{1,2,3\}$};
\node[fac] (q1) at (12.8,1.1) {1};
\node[fac] (q2) at (14.4,1.1) {2};
\node[fac] (q3) at (13.6,0) {3};
\node[fac] (q4) at (15.0,0) {4};
\draw[very thick,cA] (q1) -- node[above,font=\scriptsize] {$w$} (q2);
\draw[very thick,cA] (q1) -- node[left,font=\scriptsize] {$w$} (q3);
\draw[very thick,cA] (q2) -- node[right,font=\scriptsize] {$w$} (q3);
\node[note,anchor=north] at (13.6,-0.7) {equal flows $w$ in $Q$\\ $\Rightarrow$ enumerate location sets};
\end{tikzpicture}
\caption{Structures used in Sections~\ref{sec:exact} and~\ref{sec:clique}.}\label{fig:struct}
\end{figure}
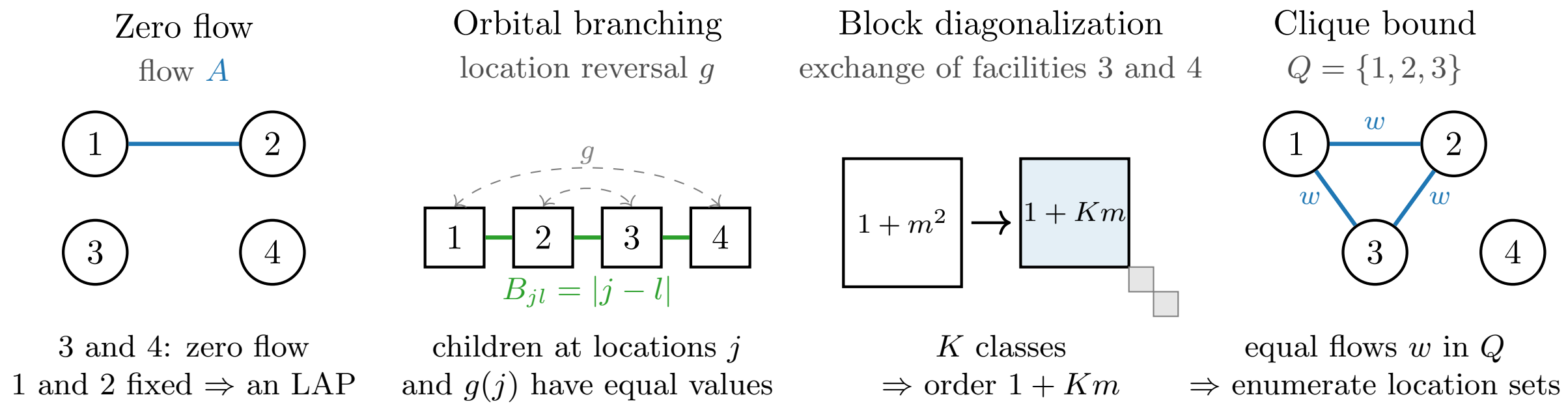

We define facilities with zero flow and the automorphisms of the data and use them to reduce the number of
nodes. None of these reductions changes the optimal value. Section~\ref{sec:num-sources} shows on which QAPLIB instances they act, and
Figure~\ref{fig:struct} shows the structures on a small example.

\begin{definition}[Zero flow]\label{def:zero}
A facility $i$ has zero flow if $A_{ik}=A_{ki}=0$ for all $k\neq i$, and a location $j$ has zero distance if
$B_{jl}=B_{lj}=0$ for all $l\neq j$.
\end{definition}

\paragraph{Facilities with zero flow.} If a facility $i$ has zero flow, its location enters the objective only through the linear term $C_{ij}+A_{ii}B_{jj}$. Such facilities  frequently appear in QAPLIB. For example, $23$ of the $32$ facilities of esc32e have zero flow. We do not branch on them~\cite{nyberg2012new,fischetti2012three}. Once all facilities with flow are fixed, or whenever $A$ or $B$ of a node is diagonal, the node is the linear assignment problem with cost matrix $C+\diag(A)\diag(B)^\top$, and \quads\ solves it exactly. Locations with zero distance become facilities with zero flow after the exchange of Section~\ref{sec:node}, which \quads\ makes at the root if it gives more facilities with zero flow.

\begin{definition}[Automorphisms]\label{def:aut}
An automorphism of a node is a pair $(\sigma,g)$ of permutations of its free facilities and its free locations
such that
\[
A_{\sigma(i)\sigma(k)}=A_{ik},\qquad B_{g(j)g(l)}=B_{jl},\qquad C_{\sigma(i)g(j)}=C_{ij}\qquad\text{for all free } i,k,j,l,
\]
and $(\sigma(i),g(j))$ is allowed whenever $(i,j)$ is allowed. The automorphisms form a group $G$. Two free facilities
$i$ and $k$ are interchangeable if $(\tau_{ik},{\rm id})\in G$, where $\tau_{ik}$ exchanges $i$ and $k$.
\end{definition}

\paragraph{Automorphisms.} Replacing $\pi$ by $g\circ\pi\circ\sigma^{-1}$ for an automorphism $(\sigma,g)$ leaves
the objective value unchanged. For example, two facilities with zero flow, equal diagonal entries of $A$, equal rows of
$C$, and the same allowed locations are interchangeable. For $Y\in\S^N$, the matrix $(\sigma,g)\cdot Y$ is obtained by permuting the rows and the columns of $Y$ so that the index $0$ is fixed and the assignment $(i,j)$ moves to $(\sigma(i),g(j))$. A matrix $Y$ is called invariant under $G$ if $(\sigma,g)\cdot Y=Y$
for every pair of $G$.

At every node, \quads\ uses the automorphisms $({\rm id},g)$ of the root with $g(j)=j$ for every location $j$ that is assigned at the node.
These automorphisms are inexpensive to enumerate. \quads\ branches on a row or a column of $X$, called a line, and creates one child for each allowed assignment of this line, up to the automorphisms (Section~\ref{sec:branch}). If the full group $G$ of the root leaves fewer children of the root than these automorphisms,
\quads\ also computes $G$ from the data of the nodes of small depth and uses it there.

\paragraph{Orbital branching.} An automorphism with $\sigma(i)=i$ maps the child that places facility $i$ at $j$ to the child that places $i$ at $g(j)$, and both children have the same optimal value. We therefore keep one child per orbit of the locations under these automorphisms~\cite{ostrowski2011orbital}, at every node of the tree. In Figure~\ref{fig:struct}, two of the four children suffice. Branching on a location is treated in the same way, and every omitted child is recorded with an automorphism that maps it to a kept child. An assignment excluded by dual fixing is excluded together with its orbit, so that the allowed assignments remain invariant under these automorphisms.

\paragraph{Duplicate and dominated nodes.} Nodes in different subtrees can also be equivalent. Let $\Gamma$ be the group of the automorphisms of the root in which $\sigma$ only permutes interchangeable facilities. A node is a duplicate of another node if an automorphism in $\Gamma$ maps the fixed assignments of the other node to those of the node. It is dominated if an automorphism in $\Gamma$ maps the fixed assignments of a node pruned by its bound into those of the node. \quads\ prunes duplicate and dominated nodes~\cite{margot2002pruning}. Since $\pi\mapsto g\circ\pi\circ\sigma^{-1}$ preserves the objective value, a duplicate node has the same optimal value as the other node. Since $\theta$ does not increase, a dominated node contains no permutation with objective value at most $\theta$. Both rules are used together with orbital branching.

\subsection{Block diagonalization and the clique bound}\label{sec:clique}

\paragraph{Block diagonalization.} The automorphisms decompose the semidefinite constraint of the DNN relaxation of a node into smaller ones~\cite{deklerk2010exploiting}. Suppose that the free facilities are partitioned into $K$ classes $\mathcal I_1,\dots,\mathcal I_K$ of interchangeable facilities, where a class may consist of only one facility. Let $E\in\{0,1\}^{m\times K}$ have the entry $E_{is}=1$ if and only if $i\in\mathcal I_s$, and let $\bar Y\coloneqq\Phi^\top Y\Phi\in\S^{1+Km}$ with $\Phi\coloneqq\begin{pmatrix}1&0\\0&I_m\otimes E\end{pmatrix}$.

\begin{proposition}\label{prop:block}
Let $Y\in\S^N$ be invariant under the permutations of the facilities within each class, nonnegative, and satisfy the gangster constraints, the arrow equalities $Y_{pp}=Y_{0p}$, and the lifted assignment
equalities. Then $Y\succeq0$ if and only if $\bar Y\succeq0$.
\end{proposition}
\begin{proof}
For two facilities $i\neq k$ of one class, let $D\coloneqq(Y_{(i,j),(i,l)})_{j,l}$ and $O\coloneqq(Y_{(i,j),(k,l)})_{j,l}$, which do not depend on the choice of $i$ and $k$ since $Y$ is invariant. By the symmetry-adapted basis of the product of the symmetric groups of the classes~\cite{deklerk2010exploiting}, $Y$ is congruent to the direct sum of a diagonal scaling of $\bar Y$ and of copies of $D-O$, one for each class with at least two facilities. The gangster constraints and the arrow equalities give
$D=\Diag(x^{(i)})$ and $O_{jj}=0$ for all $j$, where $x^{(i)}$ is the part of $x_Y$ of the facility $i$, and the lifted
assignment equalities give $Oe=x^{(i)}$. Hence $D-O=\Diag(Oe)-O$ is a weighted graph Laplacian since $O\ge0$, and it is
positive semidefinite.
\end{proof}

Hence only the block $\bar Y$ of order $1+Km$ needs to be imposed, which we call the block form of the node. Since the relaxation of the node is invariant under these permutations, averaging an optimal solution over them gives an invariant optimal solution~\cite{deklerk2010exploiting}. \quads\ uses the block form at the nodes with a large class. On esc32c, it reduces the order of the semidefinite constraint from $962$ to $280$ at the first node with a DNN bound, where $m=31$.

\paragraph{Clique bound.} On some instances, the flows are constant on groups of facilities. A set $Q$ of facilities is called a clique of flow if the flow between any two of its facilities is the same. Suppose that $C=0$ and that $A$ is
symmetric, after exchanging $A$ with $B$ if needed. For a set $Q$ of facilities, let $K_Q\in\S^m$ have the entry one at the pairs of distinct facilities of $Q$
and zero elsewhere. For a set $S$ of locations, let $B(S)\coloneqq\sum_{j,l\in S,\,j\neq l}B_{jl}$. If
$A=\sum_\iota w_\iota K_{Q_\iota}+R+\Diag(\diag(A))$ with $w_\iota>0$ and a nonnegative matrix $R$ with zero diagonal, then the objective value
of every permutation $\pi$ is
\[
\sum_\iota w_\iota\,B(\pi(Q_\iota))+\sum_{i\neq k}R_{ik}B_{\pi(i)\pi(k)}+\sum_iA_{ii}B_{\pi(i)\pi(i)},
\]
where $\pi(Q_\iota) = \{ \pi(i) \mid i\in Q_\iota\}$.
Hence, by splitting the flow matrix as in~\cite{fischetti2012three}, $v^{\rm QAP}\ge\sum_\iota w_\iota b_{|Q_\iota|}+v_R$, where $b_s\coloneqq\min\{B(S):|S|=s\}$
and $v_R$ is the GLB of the remaining two terms. \quads\ chooses the sets $Q_\iota$ greedily. It computes each $b_s$ by an enumeration of the sets of $s$ locations with pruning and with one set per orbit of the automorphisms of $B$. We call the resulting bound the clique bound. If it equals $\bar v$, the incumbent is optimal and no tree is searched. At every node, \quads\ also minimizes each term $w_\iota B(\pi(Q_\iota))$ exactly over the permutations of the node by an enumeration of the locations of the free facilities of $Q_\iota$. It adds these minima to the GLB of the remaining terms at the node and prunes the node if the sum exceeds $\theta$. If there is a single set $Q$, $R=0$, and the diagonal of $A$ is constant
on $Q$ and zero elsewhere, the objective value depends only on $S=\pi(Q)$. The same enumeration then gives the optimal
value.

\subsection{Branching rule}\label{sec:branch}

\paragraph{Choice of the line.} A line of $X$ is a row, which corresponds to a facility, or a column, which corresponds to a location. At a node, \quads\ branches on a line of $X$ and creates
one child per orbit of allowed assignments of the line. For a line $\ell$, let $E_\ell$ be the set of its
allowed assignments, with one assignment per orbit of Section~\ref{sec:exact}. Thus $E_\ell$ indexes the children
of~$\ell$. Let
\[
n_\ell\coloneqq|E_\ell|,\qquad s_\ell\coloneqq\sum_{p\in E_\ell}\max\{h,\underline v^1_p\},
\]
where $\underline v^1_p$ is the bound of child $p$ known before
it is solved (Section~\ref{sec:screen}). By default, \quads\ selects the line that minimizes $(n_\ell,-s_\ell)$ lexicographically~\cite{brixius2001solving}, so that it prefers few children and, among
them, children with high bounds. 

Two refinements apply at the top of the tree. First, at depths zero and one, if the root has no detected automorphism and $A$ and $B$ are symmetric, \quads\ improves the QPB of the children of every candidate line by a fixed number of warm-started Frank--Wolfe steps and chooses the line whose children are closest to being pruned~\cite{brixius2001solving,anstreicher2002solving,achterberg2005branching}. Second, at the nodes of depth at most two whose relative gap and that of the root are small, \quads\ selects the line with the fewest allowed assignments $p$ at which $(x_Y)_p$ is not small, where $x_Y$ is the block of the final $Y$ of the node. 
At the root, this second refinement is used only together with the first.

\paragraph{Branching on an assignment.} A node whose lower bound is close below $\theta$ is usually pruned after few further fixed assignments. However, branching on a line creates up to $m$ children. For such nodes, \quads\ branches on a single assignment. Let $a$ be the
allowed assignment with the largest diagonal entry $Y_{aa}$ of the final $Y$ of the node. If the
lower bound $h$ of the node satisfies $0\le\theta-h<\nu\max\{1,|\theta|\}\min\{1,2Y_{aa}\}$ with a small constant $\nu>0$, and more than $80\%$ of the assignments of
the node are allowed, then \quads\ splits the node into the child with $X_a=1$ and
the exclusion child with $X_a=0$. In the exclusion child, $a$ and its images in the same line under the automorphisms of the node are excluded. No objective value of the node is lost, since an automorphism maps a permutation that uses an image of $a$ to one with the same objective value that uses $a$. 
The rule is not applied at the root or at any exclusion child. 

The exclusion child fixes no facility, so it keeps the depth and the order of the matrices of its parent. It is first tested by the
lower bound of Proposition~\ref{prop:cert} at the final multiplier of the parent, restricted to its allowed assignments. Its DNN relaxation is solved only if this test does not prune it. The rules for duplicate and dominated nodes are not used in the subtree of an exclusion child.

\section{Implementation}\label{sec:impl}

We describe how \quads\ evaluates its bounds in floating-point arithmetic and how it uses several
CPU cores and GPUs in one search. This addresses the last difficulty of Section~\ref{sec:challenges}.

\subsection{Validity of the bounds in floating-point arithmetic}\label{sec:verify}

A node may be pruned only if its lower bound is valid despite rounding errors. Since
Proposition~\ref{prop:cert} gives a valid bound for every multiplier, rounding errors in ADMM affect only the
quality of the bound. Hence we run ADMM in single precision. Only the bound~\eqref{eq:phi} is evaluated in
double precision, at $\sym(Z)$. We take $t=\hat\lambda+\rho$. Here $\hat\lambda$ is a Lanczos estimate of $\lambda_{\max}(PZP)$ for
which the Cholesky factorization of $\hat\lambda I-PZP$ succeeds, and $\rho$ bounds the rounding
errors of forming $PZP$ and of the factorization~\cite{rump2006verification}. Hence $t\ge\lambda_{\max}(PZP)$. The assignment problems in~\eqref{eq:H}
and~\eqref{eq:GL2} are bounded from below by feasible dual solutions~\cite{neumaier2004safe}.

A node is pruned if $h>\theta+\epsilon\max\{1,|h|\}$, where the margin $\epsilon>0$ covers the remaining rounding
errors. Hence the incumbent is optimal when the search ends (Section~\ref{sec:search}), despite the rounding errors.

\subsection{Parallel search on CPU cores and GPUs}\label{sec:parallel}

\quads\ processes nodes in parallel on CPU cores and GPUs, whose workers share one pool of open nodes.

\paragraph{CPU workers.} \quads\ runs at most one worker per core. Each worker repeatedly removes a node
with the smallest lower bound from the pool, processes it on one core as in Algorithm~\ref{alg:bnb}, and puts its
children into the pool. The workers also share the incumbent and the statistics that control the computation of the DNN bound and dual
fixing. Since a node is pruned only by a valid bound and every node is either pruned or branched,
the optimal value and the proof of optimality do not depend on the order in which the workers process the nodes.
This order affects the number of nodes and the running time, which therefore vary slightly between runs.
The trees of the QAP are narrow at the top, where the nodes are most expensive. When fewer nodes than workers are being
processed, the idle cores are therefore used as additional threads in the matrix computations of these nodes. In particular, the root is processed by a single worker with these threads.

\paragraph{GPU workers.} \quads\ runs one worker per GPU, which also takes its nodes from the same pool. Since one node is too small to use a GPU efficiently, the worker
processes several nodes at a time and performs their ADMM iterations together. In particular, the eigenvalue
decompositions of all matrices of the same order are computed by one call of a batched eigensolver. The iterations run
in single precision on the GPU. In parallel, CPU threads evaluate the lower bounds in double precision as in
Section~\ref{sec:verify}. At the top of the tree, the pool contains few nodes. The matrices of the root have order $N$ and are
large enough to use a GPU efficiently on large instances, and the root is then solved on the GPU alone.

\section{Numerical experiments}\label{sec:num}

We compare \quads\ with solvers that compute the optimal value on QAPLIB and report its results on hard instances. We then show the effect of the structure of the data and of dual fixing.

\paragraph{Computational environment.} All experiments are run on a workstation with an Intel Xeon w7-2595X
processor with $26$ cores, $512$ GB of memory, and two NVIDIA RTX 6000 Ada Generation GPUs with $48$ GB each. \quads\
uses Intel oneMKL and CUDA~12.8. Each run uses $1$, $4$, or $16$ of these cores, as stated in each subsection. The time limit is $3600$ seconds (one hour) in Section~\ref{sec:num-compare} and $14400$ seconds (four hours) in
Sections~\ref{sec:num-dnn} and~\ref{sec:num-hard}. Every solver receives the optimal value plus one as its
upper bound~\cite{fujii2021solving}. Hence every solver must
find an optimal permutation itself.

\paragraph{Table notation.} The tables report the wall-clock time in seconds and the number of nodes. The letters K and M stand for
$10^3$ and $10^6$. The tables also report the relative gap to the optimal value $100(v^\ast-v_{\mathrm{LB}})/v^\ast$ in percent. Here $v^\ast$ is the optimal value and $v_{\mathrm{LB}}$ is
the final lower bound of the solver. For \quads, $v_{\mathrm{LB}}$ is the smallest lower bound of the open nodes when the run stops, or the optimal value if no node is open, and every bound is valid in floating-point arithmetic (Section~\ref{sec:verify}). The root gap of a bound is the relative gap with $v_{\mathrm{LB}}$ replaced by the value of the bound at the root. An underlined time marks a run stopped by the time limit, and bold marks the fastest solver in Table~\ref{tab:main}. The last row of Table~\ref{tab:main} gives the shifted geometric mean of the times with a shift of $10$ seconds, in which a
stopped run counts with the time limit, the shifted geometric mean of the numbers of nodes with a shift of $1$, and the arithmetic mean of the relative gaps. \quads\ uses the same parameters in all runs.

\subsection{Comparison with other solvers}\label{sec:num-compare}

We compare \quads\ with two solvers whose codes are available. They are of different types, and we run them on five QAPLIB instances of similar size from different families: nug25, tai25b, bur26a, kra30a, and esc32h.
\begin{itemize}[leftmargin=*,itemsep=1pt]
\item toulbar2~\cite{sewa2026snc}\footnote{\url{https://github.com/toulbar2/toulbar2}, commit \texttt{3df41536}.}, a cost function
  network solver for weighted constraint satisfaction problems, with the formulation and the options recommended for
  the QAP by its authors.
\item Gurobi~13.0.3~\cite{gurobi} on the standard binary quadratic formulation
  $\min\{x^\top(B\otimes A)x\mid \mathcal{A}x=e,\ x\in\{0,1\}^{n^2}\}$, with its default settings and a zero optimality gap.
\end{itemize}
All solvers run on one core, and Table~\ref{tab:main} reports the results.

\begin{table}[!htbp]
\centering
\caption{Comparison with other solvers on one core.}\label{tab:main}
\small
\setlength{\tabcolsep}{3pt}
\begin{tabular}{@{}l r@{\hspace{12pt}}rrr@{\hspace{12pt}}rrr@{\hspace{12pt}}rrr@{}}
\toprule
 & & \multicolumn{3}{c@{\hspace{12pt}}}{\quads} & \multicolumn{3}{c@{\hspace{12pt}}}{toulbar2} & \multicolumn{3}{c}{Gurobi}\\
\cmidrule(l{0.3em}r{\dimexpr12pt+0.3em}){3-5}\cmidrule(l{0.3em}r{\dimexpr12pt+0.3em}){6-8}\cmidrule(l{0.3em}r{0.3em}){9-11}
Instance & $n$ & time & nodes & gap & time & nodes & gap & time & nodes & gap\\
\midrule
nug25  & 25 & \textbf{69.8} & 288 & 0 & \underline{3600.0} & 21.0M & 17.79 & \underline{3600.0} & 854K & 100 \\
tai25b & 25 & \textbf{20.4} & 22 & 0 & \underline{3600.0} & 6.15M & 54.21 & \underline{3600.0} & 446K & 100 \\
bur26a & 26 & \textbf{6.3} & 27 & 0 & \underline{3600.0} & 13.0M & 1.67 & \underline{3600.0} & 286K & 100 \\
kra30a & 30 & \textbf{77.2} & 93 & 0 & \underline{3600.0} & 6.02M & 16.31 & \underline{3600.0} & 358K & 100 \\
esc32h & 32 & \textbf{0.2} & 0 & 0 & \underline{3600.0} & 63.7M & 37.44 & \underline{3600.0} & 28.9K & 100 \\
\midrule
mean & & \textbf{22.3} & 27 & 0 & 3600.0 & 14.5M & 25.48 & 3600.0 & 257K & 100 \\
\bottomrule
\end{tabular}
\end{table}

\quads\ solves all five instances within the time limit, whereas the other solvers solve none of them.
Both toulbar2 and Gurobi process many more nodes than \quads\ but do not close the gap.

\subsection{Comparison with the code of Goux}\label{sec:num-dnn}

The recent public code of J.-P.~Goux\footnote{\url{https://github.com/jpgoux/nug30}, commit \texttt{4dbaac3} of 31 October
2026.} is a branch-and-bound solver for the QAP. At every node, it computes the QPB of
Anstreicher and Brixius~\cite{anstreicher2001new}. At the nodes near the top of the tree that this bound does not prune,
it also computes a DNN bound. Its documentation and test data are restricted to the
Nugent instances, with nug30 as its target. We therefore compare on nug20, nug22, nug24, nug25, nug27, nug28, and nug30. We
run it as released, with the configuration that its documentation reports for nug30. Since its documentation reports
runs on four cores, both codes use four cores. Table~\ref{tab:nug} reports the results and the ratios of the times and of the numbers of nodes of Goux to those of \quads.

\begin{table}[!htbp]
\centering
\caption{Comparison with the code of Goux on four cores.}\label{tab:nug}
\small
\setlength{\tabcolsep}{4pt}
\begin{tabular}{@{}l r rr rr rr@{}}
\toprule
 & & \multicolumn{2}{c}{\quads} & \multicolumn{2}{c}{Goux} & \multicolumn{2}{c}{ratio}\\
\cmidrule(l{0.3em}r{0.3em}){3-4}\cmidrule(l{0.3em}r{0.3em}){5-6}\cmidrule(l{0.3em}){7-8}
Instance & $n$ & time & nodes & time & nodes & time & nodes\\
\midrule
nug20 & 20 & 3.0 & 84 & 14.0 & 141 & \textbf{4.7} & 1.7 \\
nug22 & 22 & 4.2 & 47 & 26.0 & 690 & \textbf{6.2} & 14.7 \\
nug24 & 24 & 9.6 & 190 & 68.1 & 1526 & \textbf{7.1} & 8.0 \\
nug25 & 25 & 22.6 & 286 & 205.2 & 5886 & \textbf{9.1} & 20.6 \\
nug27 & 27 & 26.1 & 197 & 315.4 & 63.3K & \textbf{12.1} & 321.3 \\
nug28 & 28 & 139.1 & 1006 & 845.5 & 44.7K & \textbf{6.1} & 44.4 \\
nug30 & 30 & 417.8 & 2841 & 3932.4 & 421K & \textbf{9.4} & 148.2 \\
\bottomrule
\end{tabular}
\end{table}

Both codes solve all seven instances, and \quads\ is several times faster on every one of them. From nug22 on, the
ratio of the numbers of nodes exceeds that of the times. Hence a node of \quads\ costs more than a node of the code of
Goux, and the advantage of \quads\ comes from its smaller tree.

\subsection{Instances solved at the root}\label{sec:num-root}
Tables~\ref{tab:root} and~\ref{tab:rootbound} list several QAPLIB instances that \quads\ solves at the root, with $16$ cores and two GPUs. Table~\ref{tab:root} lists the instances closed by the structure of the data. In this table, the entry ``cliques'' denotes the clique bound with several sets $Q_\iota$, and ``single clique'' denotes the enumeration over the sets $S$ for a single set $Q$ (Section~\ref{sec:clique}). Table~\ref{tab:rootbound} lists the instances closed by a bound at the root. In this table, RLT1 denotes the level-1 RLT bound of Section~\ref{sec:rlt1}. Most of the methods that report the shortest times on these instances have no public code. Both tables therefore give the shortest published time for each instance and the hardware used for it. These times come from different hardware and indicate only the order of magnitude. Times below $0.1$ seconds are shown as $0.1$. A dash means that no published time is known to us.

\begin{table}[!htbp]
\centering
\caption{QAPLIB instances closed by the structure of the data.}\label{tab:root}
\small
\setlength{\tabcolsep}{4pt}
\begin{tabular}{@{}l r l rr rllr@{}}
\toprule
 & & & \multicolumn{2}{c}{\quads} & \multicolumn{4}{c}{shortest published}\\
\cmidrule(l{0.3em}r{0.3em}){4-5}\cmidrule(l{0.3em}){6-9}
Instance & $n$ & structure & time & nodes & time & hardware & ref. & year\\
\midrule
esc16b & 16 & cliques & 0.1 & 0 & 2.1 & 8 cores & \cite{sewa2026snc} & 2026\\
esc16c & 16 & cliques & 0.1 & 0 & 73.1 & 8 cores & \cite{sewa2026snc} & 2026\\
esc16d & 16 & cliques & 0.1 & 0 & 0.5 & 4 cores & \cite{fischetti2012three} & 2012\\
esc16f & 16 & zero flow & 0.1 & 0 & 0.1 & 4 cores & \cite{fischetti2012three} & 2012\\
esc16h & 16 & cliques & 0.1 & 0 & 0.1 & 8 cores & \cite{sewa2026snc} & 2026\\
esc16i & 16 & cliques & 0.1 & 0 & 0.2 & 4 cores & \cite{fischetti2012three} & 2012\\
esc32e & 32 & cliques & 0.1 & 0 & 0.1 & 4 cores & \cite{fischetti2012three} & 2012\\
esc32g & 32 & cliques & 0.1 & 0 & 0.1 & 4 cores & \cite{fischetti2012three} & 2012\\
esc32h & 32 & cliques & 0.2 & 0 & 7799.9 & 4 cores & \cite{fischetti2012three} & 2012\\
esc64a & 64 & cliques & 0.1 & 0 & 83.5 & 4 cores & \cite{fischetti2012three} & 2012\\
tai64c & 64 & single clique & 1.7 & 0 & 15.0 & 1 core & \cite{drezner2015exact} & 2015\\
esc128 & 128 & cliques & 0.1 & 0 & 9.6 & 4 cores & \cite{fischetti2012three} & 2012\\
\bottomrule
\end{tabular}
\end{table}

\begin{table}[!htbp]
\centering
\caption{QAPLIB instances closed by a bound at the root.}\label{tab:rootbound}
\small
\setlength{\tabcolsep}{4pt}
\begin{tabular}{@{}l r l rr rllr@{}}
\toprule
 & & & \multicolumn{2}{c}{\quads} & \multicolumn{4}{c}{shortest published}\\
\cmidrule(l{0.3em}r{0.3em}){4-5}\cmidrule(l{0.3em}){6-9}
Instance & $n$ & bound & time & nodes & time & hardware & ref. & year\\
\midrule
chr25a & 25 & DNN & 7.3 & 1 & 7.7 & 1 core & \cite{sewa2026snc} & 2026\\
\midrule
lipa50a & 50 & RLT1 & 4.7 & 1 & 57.6 & 1 core & \cite{sewa2026snc} & 2026\\
lipa60a & 60 & RLT1 & 6.3 & 1 & 22.1 & 1 core & \cite{sewa2026snc} & 2026\\
lipa70a & 70 & RLT1 & 16.4 & 1 & 182.7 & 1 core & \cite{sewa2026snc} & 2026\\
lipa80a & 80 & RLT1 & 20.2 & 1 & -- & -- & -- & --\\
lipa90a & 90 & RLT1 & 46.0 & 1 & 106.9 & 1 core & \cite{sewa2026snc} & 2026\\
\midrule
lipa50b & 50 & GLB & 4.0 & 1 & 8.2 & 1 core & \cite{sewa2026snc} & 2026\\
lipa60b & 60 & GLB & 4.9 & 1 & 13.1 & 1 core & \cite{sewa2026snc} & 2026\\
lipa70b & 70 & GLB & 7.0 & 1 & 26.0 & 1 core & \cite{sewa2026snc} & 2026\\
lipa80b & 80 & GLB & 9.5 & 1 & 49.9 & 1 core & \cite{sewa2026snc} & 2026\\
lipa90b & 90 & GLB & 14.5 & 1 & 87.2 & 1 core & \cite{sewa2026snc} & 2026\\
\bottomrule
\end{tabular}
\end{table}

The structure of the data closes the gap on instances with up to $128$ facilities, most of them by the clique bound. On tai64c, $13$ of the $64$ facilities form a single set $Q$ as in
Section~\ref{sec:clique}, and the $512$ automorphisms of $B$ are those of an $8\times8$ torus. On esc64a, a set of $11$
facilities and a set of $3$ facilities give a clique bound equal to the optimal value. At the root, the GLB closes the gap on the lipa b instances, the level-1 RLT bound closes it on the lipa a instances, and the DNN bound closes it on chr25a. \quads\ is at least as fast as the shortest published time on all these instances.

\subsection{Instances solved by search}\label{sec:num-hard}

We run \quads\ alone with $16$ cores and both GPUs on the QAPLIB instances whose shortest published time exceeds five
minutes. Table~\ref{tab:hard}
also gives the shortest published time for each instance and the hardware used for it. This time is the wall-clock time, or the time summed over all cores where marked $^*$.

\begin{table}[!htbp]
\centering
\caption{QAPLIB instances solved by search.}\label{tab:hard}
\small
\setlength{\tabcolsep}{4pt}
\begin{tabular}{@{}l r rr rllr@{}}
\toprule
 & & \multicolumn{2}{c}{\quads} & \multicolumn{4}{c}{shortest published}\\
\cmidrule(l{0.3em}r{0.3em}){3-4}\cmidrule(l{0.3em}){5-8}
Instance & $n$ & time & nodes & time & hardware & ref. & year\\
\midrule
bur26a & 26 & 8 & 30 & 1052 & 1 core & \cite{sewa2026snc} & 2026\\
bur26b & 26 & 9 & 28 & 1169 & 1 core & \cite{sewa2026snc} & 2026\\
esc32a & 32 & 27 & 477 & 117850 & 4 cores & \cite{nyberg2013improved} & 2013\\
esc32b & 32 & 19 & 361 & 210045 & 4 cores & \cite{nyberg2013improved} & 2013\\
esc32c & 32 & 60 & 2870 & 1156 & 4 cores & \cite{fischetti2012three} & 2012\\
esc32d & 32 & 8 & 13 & 473 & 4 cores & \cite{fischetti2012three} & 2012\\
kra30a & 30 & 22 & 80 & 940 & 1 GPU & \cite{goncalves2017gpu} & 2017\\
kra30b & 30 & 38 & 1396 & 16552 & 1 GPU & \cite{goncalves2017gpu} & 2017\\
kra32 & 32 & 28 & 121 & 1059264 & $\approx$576 cores & \cite{anstreicher2002solving} & 2002\\
nug18 & 18 & 6 & 36 & 418 & 23 cores & \cite{liao2016branch} & 2016\\
nug21 & 21 & 7 & 61 & 346 & 8 cores & \cite{fujii2021solving} & 2021\\
nug25 & 25 & 13 & 270 & 491 & 1 GPU & \cite{goncalves2017gpu} & 2017\\
nug27 & 27 & 14 & 196 & 360 & 1 GPU & \cite{goncalves2017gpu} & 2017\\
nug28 & 28 & 28 & 1068 & 2084 & 1 GPU & \cite{goncalves2017gpu} & 2017\\
nug30 & 30 & 67 & 2788 & 3140 & 1728 cores & \cite{fujii2021solving} & 2021\\
rou20 & 20 & 14 & 1373 & 3419 & 23 cores & \cite{liao2016branch} & 2016\\
ste36a & 36 & 38 & 32 & 669600 & 1 core & \cite{brixius2004steinberg} & 2004\\
ste36b & 36 & 24 & 11 & 5184000\rlap{$^*$} & 22 cores & \cite{nystrom2001solving} & 2001\\
ste36c & 36 & 38 & 21 & 17280000\rlap{$^*$} & 22 cores & \cite{nystrom2001solving} & 2001\\
tai17a & 17 & 6 & 202 & 360 & 32 processors & \cite{brungger1998zram} & 1998\\
tai25a & 25 & 4344 & 218K & 109011 & 1 GPU & \cite{goncalves2017gpu} & 2017\\
tai30b & 30 & 27 & 49 & 448 & 1 GPU & \cite{goncalves2017gpu} & 2017\\
tai35b & 35 & 8129 & 180K & 249000 & 1728 cores & \cite{fujii2021solving} & 2021\\
tai40b & 40 & 272 & 161 & 105000 & 1728 cores & \cite{fujii2021solving} & 2021\\
tho30 & 30 & 306 & 15.1K & 17219 & 1 GPU & \cite{goncalves2017gpu} & 2017\\
\bottomrule
\end{tabular}
\end{table}

\quads\ solves every instance of Table~\ref{tab:hard} more than ten times faster than the shortest published time.
It uses $16$ cores and two GPUs, whereas the shortest published times of nug30, tai35b, and tai40b were obtained on
$1728$ cores. \quads\ solves nug30 in about one minute and tai40b in less than five minutes, whereas the times on
$1728$ cores are close to one hour and more than one day~\cite{fujii2021solving}. Most instances need fewer than $2000$
nodes, and only tai25a and tai35b take more than one hour. Appendix~\ref{app:all} lists the results of \quads\ on all QAPLIB instances that it solves.

\subsection{Structure of the QAPLIB families}\label{sec:num-sources}

The running times of \quads\ in Tables~\ref{tab:main}, \ref{tab:nug}, and~\ref{tab:hard} vary by orders of magnitude between instances of
similar sizes, since each family of QAPLIB has its own structure and each component of \quads\ exploits a different
one. Table~\ref{tab:struct} lists one instance of several families, run on $16$ cores and two GPUs, and gives the
structure of its data. The next four columns concern the techniques of Section~\ref{sec:reduce}. A check mark means that the run uses the technique of the column.
\begin{itemize}[leftmargin=*,itemsep=0pt]
\item Zero flow: the treatment of facilities with zero flow and locations with zero distance (Section~\ref{sec:exact}).
\item Orbits: orbital branching and the pruning of duplicate and dominated nodes (Section~\ref{sec:exact}).
\item Cliques: the clique bound at the nodes (Section~\ref{sec:clique}).
\item Cuts: the cuts at the root (Section~\ref{sec:cuts}).
\end{itemize}
The level-1 RLT bound is used only at the root. Several instances closed at the root are listed in
Tables~\ref{tab:root} and~\ref{tab:rootbound}. The last columns give the root gaps of the GLB, the
QPB, and the DNN bound, where a root gap above $100\%$ means a negative bound. A dash marks a
bound that is not computed, such as the QPB of bur26a, whose matrices $A$ and $B$ are both asymmetric.

\begin{table}[!htbp]
\centering
\caption{Techniques of \quads\ on one instance of several families.}\label{tab:struct}
\small
\setlength{\tabcolsep}{4pt}
\begin{tabular}{@{}l l cccc rrr r@{}}
\toprule
 & & & & & & \multicolumn{3}{c}{root gap} & \\
\cmidrule(lr){7-9}
Instance & structure & zero flow & orbits & cliques & cuts & GLB & QPB & DNN & nodes\\
\midrule
esc32a & sparse flow & \checkmark & \checkmark & \checkmark & & 73.08 & 215.8 & -- & 477 \\
\midrule
ste36a & wiring, sparse flow & \checkmark & \checkmark & & & 25.22 & 211.8 & 2.99 & 32 \\
\midrule
kra32 & hospital layout & \checkmark & & & & 24.02 & 29.39 & 3.96 & 121 \\
\midrule
bur26a & keyboard layout & & \checkmark & & & 2.05 & -- & 0.14 & 30 \\
nug28 & grid distances & & \checkmark & & & 26.71 & 12.66 & 3.26 & 1068 \\
kra30b & hospital layout & & \checkmark & & & 24.45 & 26.70 & 4.41 & 1396 \\
tho30 & grid distances & & \checkmark & & & 39.59 & 17.13 & 4.85 & 15.1K \\
\midrule
tai30b & asymmetric random & & & & \checkmark & 88.01 & 45.57 & 0.03 & 49 \\
\bottomrule
\end{tabular}
\end{table}

These techniques speed up \quads\ on the instances whose data have the corresponding structure. The treatment of zero flow and orbital branching reduce the number of nodes, and the cuts reduce the root gap.

\subsection{Dual fixing and the DNN bound}\label{sec:num-dualfix}

\begin{figure}[!htbp]
\centering
\includegraphics[width=0.85\linewidth]{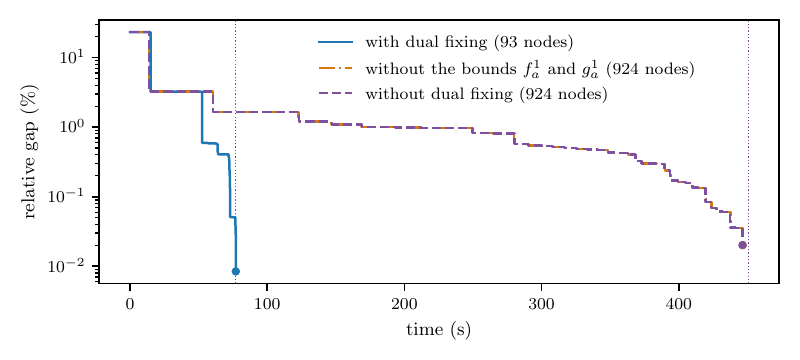}
\caption{Relative gap on kra30a with and without dual fixing.}\label{fig:dualfix}
\end{figure}

\begin{figure}[!htbp]
\centering
\includegraphics[width=\linewidth]{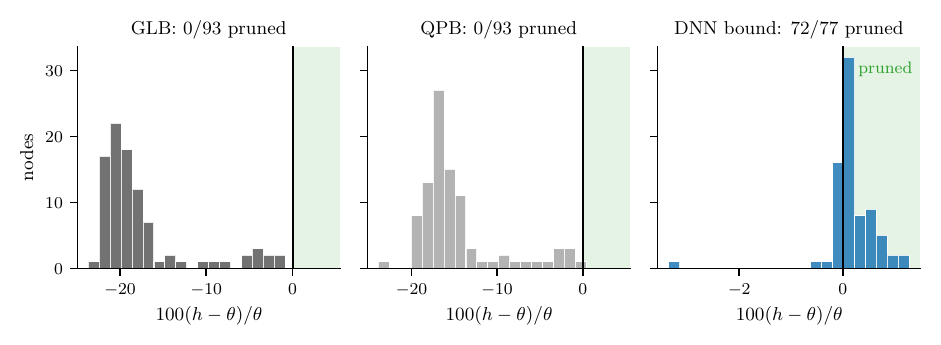}
\caption{Bounds of the nodes of kra30a.}\label{fig:bounds}
\end{figure}

Figures~\ref{fig:dualfix} and~\ref{fig:bounds} show the effect of dual fixing and of the DNN bound on kra30a, with one core for each run. Figure~\ref{fig:dualfix} compares \quads\ with two variants by the relative gap over time. The first variant omits the bounds $f^1_a$ and $g^1_a$ of Section~\ref{sec:dualfix}, and the second omits dual fixing entirely. Figure~\ref{fig:bounds} shows the distribution of $100(h-\theta)/\theta$ over the $93$ nodes processed by \quads, with the DNN bound at the $77$ nodes at which it is computed, where $h$ is the GLB, the QPB, or the DNN bound of a node. A node to the right of the line is pruned by the bound.

Dual fixing reduces the number of nodes by a factor of about ten and the running time by a factor of about six. The DNN bound prunes almost every node at which it is solved, whereas the GLB and the QPB prune none and lie far below the threshold. Hence the DNN bound, rather than the cheap bounds, drives the search and dual fixing on this instance.

\section{Conclusions}\label{sec:conclusion}

In this paper, we have proposed \quads, a branch-and-bound solver for the QAP based on DNN bounds computed by a
two-block ADMM. Since a valid lower bound is available from any multiplier of the ADMM, each node can stop early. The same multiplier warm starts the children and excludes assignments from the subtree of the node by dual fixing. The structure
of the data further reduces the search tree through the treatment of facilities with zero flow, orbital branching,
and the clique bound, and block diagonalization reduces the cost of a node. Aggregated RLT cuts at the root strengthen the
bounds of the root and of its descendants. On QAPLIB, \quads\ is the fastest of the compared solvers on one core and is faster than the code of Goux on four cores. In parallel on CPU cores and GPUs, it solves tai40b in less than five minutes. On every QAPLIB instance that \quads\ solves and whose shortest published time exceeds five minutes, it is more than ten times faster than that time. The clique bound and the treatment of zero flow prove optimality at the root for several {\tt esc} instances with up to $128$ facilities. Separating the cuts below the root and branching on classes of interchangeable facilities are left for future work.


\clearpage
\appendix
\section{Results on all solved QAPLIB instances}\label{app:all}

Table~\ref{tab:all} lists the time and the number of nodes of \quads\ on every QAPLIB instance that it solves, with $16$
CPU cores and two GPUs.

\begin{table}[!htbp]
\centering
\caption{\quads\ on all solved QAPLIB instances with $16$ cores and two GPUs.}\label{tab:all}
\footnotesize
\setlength{\tabcolsep}{3pt}
\begin{tabular}{@{}lrrr@{}p{24pt}@{}lrrr@{}p{24pt}@{}lrrr@{}}
\toprule
Instance & $n$ & time & nodes & & Instance & $n$ & time & nodes & & Instance & $n$ & time & nodes\\
\cmidrule{1-4}\cmidrule{6-9}\cmidrule{11-14}
bur26a & 26 & 8.0 & 30 & & esc32c & 32 & 59.8 & 2870 & & nug16b & 16 & 4.0 & 4\\
bur26b & 26 & 8.6 & 28 & & esc32d & 32 & 7.8 & 13 & & nug17 & 17 & 5.1 & 32\\
bur26c & 26 & 8.2 & 47 & & esc32e & 32 & 0.1 & 0 & & nug18 & 18 & 6.0 & 36\\
bur26d & 26 & 9.0 & 48 & & esc32g & 32 & 0.1 & 0 & & nug20 & 20 & 5.9 & 84\\
bur26e & 26 & 13.8 & 62 & & esc32h & 32 & 0.2 & 0 & & nug21 & 21 & 7.3 & 61\\
bur26f & 26 & 7.3 & 51 & & esc64a & 64 & 0.1 & 0 & & nug22 & 22 & 7.1 & 46\\
bur26g & 26 & 7.2 & 31 & & esc128 & 128 & 0.1 & 0 & & nug24 & 24 & 9.8 & 190\\
bur26h & 26 & 7.1 & 21 & & had12 & 12 & 3.4 & 1 & & nug25 & 25 & 12.8 & 270\\
chr12a & 12 & 3.4 & 1 & & had14 & 14 & 3.5 & 1 & & nug27 & 27 & 13.5 & 196\\
chr12b & 12 & 3.4 & 1 & & had16 & 16 & 3.5 & 1 & & nug28 & 28 & 28.3 & 1068\\
chr12c & 12 & 3.4 & 1 & & had18 & 18 & 4.2 & 1 & & nug30 & 30 & 67.2 & 2788\\
chr15a & 15 & 3.6 & 1 & & had20 & 20 & 4.3 & 1 & & rou12 & 12 & 3.9 & 1\\
chr15b & 15 & 3.5 & 1 & & kra30a & 30 & 22.2 & 80 & & rou15 & 15 & 5.1 & 27\\
chr15c & 15 & 3.5 & 1 & & kra30b & 30 & 38.1 & 1396 & & rou20 & 20 & 14.4 & 1373\\
chr18a & 18 & 4.4 & 8 & & kra32 & 32 & 28.3 & 121 & & scr12 & 12 & 3.5 & 1\\
chr18b & 18 & 3.3 & 1 & & lipa20a & 20 & 3.6 & 1 & & scr15 & 15 & 3.8 & 1\\
chr20a & 20 & 4.7 & 1 & & lipa20b & 20 & 3.4 & 1 & & scr20 & 20 & 7.0 & 26\\
chr20b & 20 & 7.4 & 1 & & lipa30a & 30 & 5.8 & 1 & & ste36a & 36 & 38.1 & 32\\
chr20c & 20 & 4.5 & 1 & & lipa30b & 30 & 3.3 & 1 & & ste36b & 36 & 23.8 & 11\\
chr22a & 22 & 4.7 & 1 & & lipa40a & 40 & 9.4 & 1 & & ste36c & 36 & 38.1 & 21\\
chr22b & 22 & 5.8 & 1 & & lipa40b & 40 & 3.6 & 1 & & tai12a & 12 & 3.4 & 1\\
chr25a & 25 & 7.3 & 1 & & lipa50a & 50 & 4.7 & 1 & & tai12b & 12 & 4.0 & 8\\
els19 & 19 & 4.3 & 32 & & lipa50b & 50 & 4.0 & 1 & & tai15a & 15 & 4.9 & 68\\
esc16a & 16 & 4.1 & 5 & & lipa60a & 60 & 6.3 & 1 & & tai15b & 15 & 4.1 & 11\\
esc16b & 16 & 0.1 & 0 & & lipa60b & 60 & 4.9 & 1 & & tai17a & 17 & 6.1 & 202\\
esc16c & 16 & 0.1 & 0 & & lipa70a & 70 & 16.4 & 1 & & tai20a & 20 & 12.5 & 1080\\
esc16d & 16 & 0.1 & 0 & & lipa70b & 70 & 7.0 & 1 & & tai20b & 20 & 12.7 & 15\\
esc16e & 16 & 3.5 & 2 & & lipa80a & 80 & 20.2 & 1 & & tai25a & 25 & 4344.0 & 218K\\
esc16f & 16 & 0.1 & 0 & & lipa80b & 80 & 9.5 & 1 & & tai25b & 25 & 15.3 & 22\\
esc16g & 16 & 3.4 & 2 & & lipa90a & 90 & 46.0 & 1 & & tai30b & 30 & 26.7 & 49\\
esc16h & 16 & 0.1 & 0 & & lipa90b & 90 & 14.5 & 1 & & tai35b & 35 & 8129.0 & 180K\\
esc16i & 16 & 0.1 & 0 & & nug12 & 12 & 3.7 & 5 & & tai40b & 40 & 271.5 & 161\\
esc16j & 16 & 3.4 & 1 & & nug14 & 14 & 4.1 & 14 & & tai64c & 64 & 1.7 & 0\\
esc32a & 32 & 26.6 & 477 & & nug15 & 15 & 3.9 & 7 & & tho30 & 30 & 306.4 & 15.1K\\
esc32b & 32 & 19.4 & 361 & & nug16a & 16 & 4.6 & 19 & & & & &\\
\bottomrule
\end{tabular}
\end{table}


\begin{thebibliography}{99}\small

\bibitem{achterberg2005branching} T.~Achterberg, T.~Koch, and A.~Martin. Branching rules revisited. \emph{Oper. Res. Lett.}, 33 (2005), pp.~42--54.


\bibitem{adams1994improved} W.~P. Adams and T.~A. Johnson. Improved linear programming-based lower bounds for the quadratic assignment problem. In \emph{Quadratic Assignment and Related Problems}, DIMACS Ser. Discrete Math. Theoret. Comput. Sci. 16, AMS, 1994, pp.~43--75.

\bibitem{adams2007level2} W.~P. Adams, M.~Guignard, P.~M. Hahn, and W.~L. Hightower. A level-2 reformulation-linearization technique bound for the quadratic assignment problem. \emph{European J. Oper. Res.}, 180 (2007), pp.~983--996.


\bibitem{anstreicher2001new} K.~M. Anstreicher and N.~W. Brixius. A new bound for the quadratic assignment problem based on convex quadratic programming. \emph{Math. Program.}, 89 (2001), pp.~341--357.

\bibitem{anstreicher2002solving} K.~M. Anstreicher, N.~W. Brixius, J.-P. Goux, and J.~Linderoth. Solving large quadratic assignment problems on computational grids. \emph{Math. Program.}, 91 (2002), pp.~563--588.



\bibitem{brixius2001solving} K.~M. Anstreicher and N.~W. Brixius. Solving quadratic assignment problems using convex quadratic programming relaxations. \emph{Optim. Methods Softw.}, 16 (2001), pp.~49--68.

\bibitem{brixius2004steinberg} N.~W. Brixius and K.~M. Anstreicher. The Steinberg wiring problem. In \emph{The Sharpest Cut}, MPS/SIAM Ser. Optim. 4, SIAM, Philadelphia, 2004, pp.~293--307.

\bibitem{brungger1998zram} A.~Br\"ungger, A.~Marzetta, J.~Clausen, and M.~Perregaard. Solving large-scale QAP problems in parallel with the search library ZRAM. \emph{J. Parallel Distrib. Comput.}, 50 (1998), pp.~157--169.

\bibitem{burkard1997qaplib} R.~E. Burkard, S.~E. Karisch, and F.~Rendl. QAPLIB -- a quadratic assignment problem library. \emph{J. Global Optim.}, 10 (1997), pp.~391--403.

\bibitem{burkard2012assignment} R.~E. Burkard, M.~Dell'Amico, and S.~Martello. \emph{Assignment Problems}, revised reprint. SIAM, Philadelphia, 2012.

\bibitem{deklerk2010exploiting} E.~de Klerk and R.~Sotirov. Exploiting group symmetry in semidefinite programming relaxations of the quadratic assignment problem. \emph{Math. Program.}, 122 (2010), pp.~225--246.

\bibitem{demeijer2023partitioning} F.~de Meijer, R.~Sotirov, A.~Wiegele, and S.~Zhao. Partitioning through projections: strong SDP bounds for large graph partition problems. \emph{Comput. Oper. Res.}, 151 (2023), Art.~106088.


\bibitem{drezner2015exact} Z.~Drezner, A.~Misevi\v{c}ius, and G.~Palubeckis. Exact algorithms for the solution of the grey pattern quadratic assignment problem. \emph{Math. Methods Oper. Res.}, 82 (2015), pp.~85--105.

\bibitem{faye2005cutting} A.~Faye and F.~Roupin. A cutting planes algorithm based upon a semidefinite relaxation for the quadratic assignment problem. In \emph{Algorithms -- ESA 2005}, LNCS 3669, Springer, 2005, pp.~850--861.

\bibitem{fazel2013hankel} M.~Fazel, T.~K. Pong, D.~Sun, and P.~Tseng. Hankel matrix rank minimization with applications to system identification and realization. \emph{SIAM J. Matrix Anal. Appl.}, 34 (2013), pp.~946--977.

\bibitem{fischetti2012three} M.~Fischetti, M.~Monaci, and D.~Salvagnin. Three ideas for the quadratic assignment problem. \emph{Oper. Res.}, 60 (2012), pp.~954--964.

\bibitem{frank1956algorithm} M.~Frank and P.~Wolfe. An algorithm for quadratic programming. \emph{Naval Res. Logist. Quart.}, 3 (1956), pp.~95--110.

\bibitem{fujii2021solving} K.~Fujii, N.~Ito, S.~Kim, M.~Kojima, Y.~Shinano, and K.-C. Toh. Solving challenging large scale QAPs. arXiv:2101.09629, 2021.

\bibitem{gilmore1962optimal} P.~C. Gilmore. Optimal and suboptimal algorithms for the quadratic assignment problem. \emph{J. Soc. Indust. Appl. Math.}, 10 (1962), pp.~305--313.

\bibitem{glowinski1975approximation} R.~Glowinski and A.~Marroco. Sur l'approximation, par \'el\'ements finis d'ordre un, et la r\'esolution, par p\'enalisation-dualit\'e, d'une classe de probl\`emes de Dirichlet non lin\'eaires. \emph{RAIRO Anal. Num\'er.}, 9 (1975), pp.~41--76.

\bibitem{goncalves2017gpu} A.~D. Gon\c{c}alves, A.~A. Pessoa, C.~Bentes, R.~Farias, and L.~M. de~A. Drummond. A graphics processing unit algorithm to solve the quadratic assignment problem using level-2 reformulation-linearization technique. \emph{INFORMS J. Comput.}, 29 (2017), pp.~676--687.

\bibitem{graham2022restricted} N.~Graham, H.~Hu, J.~Im, X.~Li, and H.~Wolkowicz. A restricted dual Peaceman--Rachford splitting method for a strengthened DNN relaxation for QAP. \emph{INFORMS J. Comput.}, 34 (2022), pp.~2125--2143.

\bibitem{gurobi} Gurobi Optimization, LLC. \emph{Gurobi Optimizer}, version 13.0.3, 2026. \url{https://support.gurobi.com/hc/en-us/articles/50235313254161-Gurobi-13-0-3-Released}.


\bibitem{hahn1998lower} P.~Hahn and T.~Grant. Lower bounds for the quadratic assignment problem based upon a dual formulation. \emph{Oper. Res.}, 46 (1998), pp.~912--922.

\bibitem{helmberg2000fixing} C.~Helmberg. Fixing variables in semidefinite relaxations. \emph{SIAM J. Matrix Anal. Appl.}, 21 (2000), pp.~952--969.
\bibitem{hildreth1957quadratic} C.~Hildreth. A quadratic programming procedure. \emph{Naval Res. Logist. Quart.}, 4 (1957), pp.~79--85.

\bibitem{hou2025rinnal} D.~Hou, T.~Tang, and K.-C. Toh. A low-rank augmented Lagrangian method for doubly nonnegative relaxations of mixed-binary quadratic programs. \emph{Oper. Res.}, 74 (2026), pp.~1626--1647.

\bibitem{hou2025rinnalplus} D.~Hou, T.~Tang, and K.-C. Toh. RiNNAL$+$: a Riemannian ALM solver for SDP-RLT relaxations of mixed-binary quadratic programs. \emph{Math. Program. Comput.}, 18 (2026), pp.~679--735.

\bibitem{hou2026scsdp} D.~Hou, T.~P.~D. Nguyen, K.-C. Toh, and G.~Wang. Sparsity-cone SDP relaxations and applications to variable fixing for sparse quadratic programs. arXiv:2606.22894, 2026.


\bibitem{hu2021linearization} H.~Hu and R.~Sotirov. The linearization problem of a binary quadratic problem and its applications. \emph{Ann. Oper. Res.}, 307 (2021), pp.~229--249.

\bibitem{hu2026exact} H.~Hu and M.~Xu. Exact convex reformulations of binary quadratic programs from facially reduced semidefinite relaxations. arXiv:2610.11154, 2026.


\bibitem{junger2001box} M.~J\"unger and V.~Kaibel. Box-inequalities for quadratic assignment polytopes. \emph{Math. Program.}, 91 (2001), pp.~175--197.


\bibitem{kim2016lagrangian} S.~Kim, M.~Kojima, and K.-C. Toh. A Lagrangian-DNN relaxation: a fast method for computing tight lower bounds for a class of quadratic optimization problems. \emph{Math. Program.}, 156 (2016), pp.~161--187.
\bibitem{koopmans1957assignment} T.~C. Koopmans and M.~Beckmann. Assignment problems and the location of economic activities. \emph{Econometrica}, 25 (1957), pp.~53--76.

\bibitem{krislock2014improved} N.~Krislock, J.~Malick, and F.~Roupin. Improved semidefinite bounding procedure for solving Max-Cut problems to optimality. \emph{Math. Program.}, 143 (2014), pp.~61--86.

\bibitem{kuhn1955hungarian} H.~W. Kuhn. The Hungarian method for the assignment problem. \emph{Naval Res. Logist. Quart.}, 2 (1955), pp.~83--97.

\bibitem{lawler1963quadratic} E.~L. Lawler. The quadratic assignment problem. \emph{Management Sci.}, 9 (1963), pp.~586--599.

\bibitem{li2020lowrank} Y.~Li, H.~Liu, Z.~Wen, and Y.-x.~Yuan. Low-rank matrix iteration using polynomial-filtered subspace extraction. \emph{SIAM J. Sci. Comput.}, 42 (2020), pp.~A1686--A1713.

\bibitem{liao2016branch} Z.~Liao. Branch and bound via the alternating direction method of multipliers for the quadratic assignment problem. Master's research paper, University of Waterloo, 2016.
\bibitem{loiola2007survey} E.~M. Loiola, N.~M.~M. de~Abreu, P.~O. Boaventura-Netto, P.~Hahn, and T.~Querido. A survey for the quadratic assignment problem. \emph{European J. Oper. Res.}, 176 (2007), pp.~657--690.

\bibitem{margot2002pruning} F.~Margot. Pruning by isomorphism in branch-and-cut. \emph{Math. Program.}, 94 (2002), pp.~71--90.

\bibitem{mittelmann2026qaplib} H.~D. Mittelmann. Improved QAPLIB lower bounds using BBCPOP and NewtBracket.
\url{https://plato.asu.edu/ftp/qaplib_bounds.html}, version of 7 November 2026.

\bibitem{neumaier2004safe} A.~Neumaier and O.~Shcherbina. Safe bounds in linear and mixed-integer linear programming. \emph{Math. Program.}, 99 (2004), pp.~283--296.


\bibitem{nyberg2012new} A.~Nyberg and T.~Westerlund. A new exact discrete linear reformulation of the quadratic assignment problem. \emph{European J. Oper. Res.}, 220 (2012), pp.~314--319.

\bibitem{nyberg2013improved} A.~Nyberg, T.~Westerlund, and A.~Lundell. Improved discrete reformulations for the quadratic assignment problem. In \emph{CPAIOR 2013}, LNCS 7874, Springer, 2013, pp.~193--203.

\bibitem{nystrom2001solving} M.~Nystr\"{o}m. Solving certain large instances of the quadratic assignment problem: Steinberg's examples. Technical Report CaltechCSTR:2001.010, California Institute of Technology, 2001.
\bibitem{oliveira2018admm} D.~E. Oliveira, H.~Wolkowicz, and Y.~Xu. ADMM for the SDP relaxation of the QAP. \emph{Math. Program. Comput.}, 10 (2018), pp.~631--658.

\bibitem{ostrowski2011orbital} J.~Ostrowski, J.~Linderoth, F.~Rossi, and S.~Smriglio. Orbital branching. \emph{Math. Program.}, 126 (2011), pp.~147--178.


\bibitem{povh2009copositive} J.~Povh and F.~Rendl. Copositive and semidefinite relaxations of the quadratic assignment problem. \emph{Discrete Optim.}, 6 (2009), pp.~231--241.

\bibitem{rendl2010solving} F.~Rendl, G.~Rinaldi, and A.~Wiegele. Solving Max-Cut to optimality by intersecting semidefinite and polyhedral relaxations. \emph{Math. Program.}, 121 (2010), pp.~307--335.

\bibitem{rockafellar1970convex} R.~T. Rockafellar. \emph{Convex Analysis}. Princeton University Press, Princeton, 1970.

\bibitem{rump2006verification} S.~M. Rump. Verification of positive definiteness. \emph{BIT Numer. Math.}, 46 (2006), pp.~433--452.

\bibitem{sahni1976pcomplete} S.~Sahni and T.~Gonzalez. P-complete approximation problems. \emph{J. ACM}, 23 (1976), pp.~555--565.

\bibitem{sewa2026snc} G.~Sewa, D.~Allouche, S.~de~Givry, G.~Katsirelos, and T.~Schiex. Singleton node consistency for quadratic assignment problems in cost function networks. In \emph{CPAIOR 2026}, LNCS 16595, Springer, 2026, pp.~520--537.


\bibitem{sherali1999reformulation} H.~D. Sherali and W.~P. Adams. \emph{A Reformulation-Linearization Technique for Solving Discrete and Continuous Nonconvex Problems}. Kluwer, Dordrecht, 1999.

\bibitem{steinberg1961backboard} L.~Steinberg. The backboard wiring problem: a placement algorithm. \emph{SIAM Rev.}, 3 (1961), pp.~37--50.

\bibitem{taillard1991robust} \'E.~D. Taillard. Robust taboo search for the quadratic assignment problem. \emph{Parallel Comput.}, 17 (1991), pp.~443--455.

\bibitem{yang2015sdpnalplus} L.~Yang, D.~Sun, and K.-C. Toh. SDPNAL$+$: a majorized semismooth Newton-CG augmented Lagrangian method for semidefinite programming with nonnegative constraints. \emph{Math. Program. Comput.}, 7 (2015), pp.~331--366.

\bibitem{yang2026inexact} L.~Yang, H.~Lan, D.~Hou, L.~Liang, and K.-C. Toh. An inexact Halpern-accelerated preconditioned generalized proximal point algorithm for the maximal monotone inclusion problem. arXiv:2610.05742, 2026.


\bibitem{zhao1998semidefinite} Q.~Zhao, S.~E. Karisch, F.~Rendl, and H.~Wolkowicz. Semidefinite programming relaxations for the quadratic assignment problem. \emph{J. Comb. Optim.}, 2 (1998), pp.~71--109.

\end{thebibliography}
\end{document}